\documentclass[onefignum,onetabnum]{siamart171218}
\pdfoutput=1

\usepackage{subfig}
\usepackage{amssymb} %
\usepackage{mathtools} %
\usepackage{graphicx}
\usepackage{color}
\usepackage{tabularx}
\usepackage{tikz} %
\usepackage{enumitem}
\usepackage[labelfont=bf]{caption}
\usepackage{changepage}
\usepackage{afterpage} %
\usepackage{gensymb} %
\usepackage{placeins} %
\usepackage{enumitem}%
\usepackage{siunitx} %
\usepackage{soul} %
\usepackage[export]{adjustbox} %
\usepackage{tabularx}
\usepackage{wrapfig}
\usepackage{mwe}
\usepackage{dashrule}
\usepackage[export]{adjustbox}
\usepackage{lscape}

\usepackage{lipsum}
\usepackage{amsfonts}
\usepackage{graphicx}
\usepackage{epstopdf}
\usepackage{algorithmic}
\ifpdf
  \DeclareGraphicsExtensions{.eps,.pdf,.png,.jpg}
\else
  \DeclareGraphicsExtensions{.eps}
\fi

\newsiamremark{remark}{Remark}
\newsiamremark{hypothesis}{Hypothesis}
\crefname{hypothesis}{Hypothesis}{Hypotheses}
\newsiamthm{claim}{Claim}

\headers{Exponential Polynomial Block Methods}{T. Buvoli}

\title{Exponential Polynomial Block Methods \thanks{Submitted to the editors DATE.
\funding{This work was funded by the National Science Foundation, Computational Mathematics Program DMS-1216732 and DMS-2012875}}}

\author{Tommaso Buvoli\thanks{Department of Applied Mathematics, University of California, Merced, Merced CA USA. 
  (\email{tbuvoli@ucmerced.edu}).}
}

\usepackage{amsopn}

\ifpdf
\hypersetup{
  pdftitle={Exponential Polynomial Block Methods}%
  pdfauthor={T. Buvoli}
}
\fi

\newcommand{\ODEd}{ODE dataset}
\newcommand{\ODED}{ODE Dataset}

\newcommand{\ODEp}{ODE polynomial}

\newcommand{\ODEps}{ODE polynomials}

\newcommand{\ODEdp}{ODE derivative polynomial}

\newcommand{\ppe}{polynomial $\varphi$-expansion}

\newcommand{\ppes}{polynomial $\varphi$-expansions}

\newcommand{\Ape}{Adams $\varphi$-expansion}
\newcommand{\APE}{Adams $\varphi$-Expansion}
\newcommand{\Apes}{Adams $\varphi$-expansions}

	\newcommand{\ODEdcp}{ODE derivative component polynomial}
	\newcommand{\ODEDCP}{ODE Derivative Component Polynomial}

\newcommand{\expDSU}[1][s]{D_R(r,#1)}					%
\newcommand{\expDSP}[1][s]{D_N(r,#1)}					%
\definecolor{plot_black}{RGB}{10, 10, 10}
\definecolor{plot_red}{RGB}{192, 57, 43}
\definecolor{plot_orange}{RGB}{230, 126, 34}
\definecolor{plot_yellow}{RGB}{241, 196, 15}
\definecolor{plot_green}{RGB}{39, 174, 96}
\definecolor{plot_blue}{RGB}{41, 128, 185}
\definecolor{plot_violet}{RGB}{155, 89, 182}
\definecolor{plot_grey}{RGB}{149, 165, 166}

\begin{document}

\maketitle

\begin{abstract}
	In this paper we extend the polynomial time integration framework to include exponential integration for both partitioned and unpartitioned initial value problems. We then demonstrate the utility of the exponential polynomial framework by constructing a new class of parallel exponential polynomial block methods (EPBMs) based on the Legendre points. These new integrators can be constructed at arbitrary orders of accuracy and have  improved stability compared to existing exponential linear multistep methods. Moreover, if the ODE right-hand side evaluations can be parallelized efficiently, then high-order EPBMs are significantly more efficient at obtaining highly accurate solutions than exponential linear multistep methods and exponential spectral deferred correction methods of equivalent order.
\end{abstract}

\begin{keywords}
	Time-integration, polynomial interpolation, general linear methods, parallelism, high-order.
\end{keywords}

\begin{AMS}
	65L04, 65L05, 65L06, 65L07
\end{AMS}

Polynomial time integrators \cite{buvoli2018polynomial, buvoli2019constructing} are a class of parametrized methods for solving first-order systems of ordinary differential equations. The integrators are based on a new  framework that combines ideas from complex analysis, approximation theory, and general linear methods. The framework encompasses all general linear methods that compute their outputs and stages using interpolating polynomials in the complex time plane. The use of polynomials enables one to trivially satisfy order conditions and easily construct a range of implicit or explicit integrators with properties such as parallelism and high-order of accuracy. 

In order to extend the utility of the polynomial framework, we generalize it to include exponential integration. %
Exponential integrators are a general class of methods that incorporate exponential functions to provide increased stability and accuracy for solving stiff systems \cite{hochbruck2010exponentialreview}. Continuing efforts to construct and analyze exponential integrators have already produced a wide range of methods \cite{beylkin1998ELP, cox2002ETDRK4, KassamTrefethen05ETDRK4, krogstad2005IF,hochbruckostermann2005ETDRKSTIFFA, hochbruckostermann2005ETDRKSTIFFB, koikari2005rooted, hochbruck2010exponentialreview, ostermann2006general} that can provide improved efficiency compared to fully implicit and semi-implicit integrators \cite{grooms2011IMEXETDCOMP, KassamTrefethen05ETDRK4, loffeld2013comparative, montanelli2016solving}. 

Incorporating parallelism into exponential integrators remains an open question.
Typically, parallelism is applied within an exponential scheme to speed up the estimation of exponential matrix functions and their products with vectors; examples include parallel Krylov projections \cite{loffeld2014implementation} and parallel rational approximations of exponential functions \cite{sheen2000parallel, haut2015high, schreiber2019parallel, schreiber2019exponential}. To the best of our knowledge there has only been limited research in developing exponential integrators that naturally incorporate parallelism in their stages and outputs. Exponential EPIRK methods with independent stages \cite{rainwater2016new} and parallel exponential Rosenbrock methods \cite{luan2016parallel} constitute exceptions to this assessment. However, both of these approaches are limited since they require a restricted integrator formulation that only permits a limited number of stages to be parallelized. Furthermore, it is difficult to derive arbitrary-order parallel schemes of this type. 

The aim of this work is therefore twofold. First, we extend the polynomial framework to include exponential integration and introduce a general formulation that encompasses all families of exponential polynomial integrators. Second, we demonstrate the utility of the framework by presenting several method constructions for deriving both serial or parallel exponential polynomial block methods (EPBMs) with any order-of-accuracy. Unlike existing exponential integrators, the new parallel method constructions enable simultaneous computation of all their output values. %

This paper is organized as follows. In Sections \ref{sec:polynomial_introduction} and \ref{sec:exponential_introduction} we provide a brief introduction to polynomial methods and exponential integrators. In Section \ref{sec:exponential_polynomial} we extend the polynomial framework to include exponential integration and introduce polynomial block methods. In Section \ref{sec:constructing_pbm} we propose several general strategies for constructing parallel or serial polynomial block methods of any order. We also introduce a new class of polynomial block methods based on the Legendre points. Next, in Section \ref{sec:stability_accuracy}, we analyze and compare the stability regions of the new methods to existing exponential Adams-Bashforth and exponential Runge-Kutta methods. Finally, in Section \ref{sec:numerical_experiments}, we perform numerical experiments comparing EPBMs to a variety of other high-order exponential integrators.

\section{Polynomial time integrators} 
\label{sec:polynomial_introduction}

Polynomial time integrators \cite{buvoli2018polynomial, buvoli2019constructing} are a general class of parametrized methods constructed using interpolating polynomials in the complex time plane. The polynomial framework is based on the {\em \ODEp{}} and the {\em \ODEd{}}. The former describes a range of polynomial-based approximations, and the latter contains the data values for constructing these interpolants. The primary distinguishing feature of a polynomial method is that each of its stages and outputs are computed by evaluating an \ODEp{}. Broadly speaking, the polynomial framework encapsulates the subset of all general linear methods with coefficients that correspond to those of ODE polynomials.

In the subsections that follow, we briefly review the \ODEd{}, the \ODEp{}, and the notation used to describe polynomial time integrators. %

\subsection{The \ODEd{}}

The \ODEd{} is an ordered set containing all possible data values for constructing the interpolating polynomials that are used to compute a method's outputs.  At the $n$th timestep an \ODEd{} contains a method's inputs, outputs, and stage values along with their derivatives and their temporal nodes. The data is represented in the local coordinates $\tau$ where the global time is
	\begin{align}
		t = r\tau + t_n,
		\label{eq:global_time}
	\end{align}
	and $r$ is a scaling factor. An \ODEd{} of size $w$ is represented with the notation 
	\begin{align}
		D(r, t_n) =\left\{ \left(\tau_j,~ y_j,~ r f_j \right) \right\}_{j=1}^w
		\label{eq:generic_ode_dataset}
	\end{align}
	where ${y_j \approx y(t(\tau_j))}$, and ${f_j = F(t(\tau_j), y_j)}$. 	
	
\subsection{The \ODEp{}}

An \ODEp{} can be used to represent a wide variety of approximations for the Taylor series of a differential equation solution. In its most general form, 	an \ODEp{} of degree $g$ with expansion point $b$ is
			\begin{align}
				p(\tau; b) &= \sum_{j=0}^{g}  \frac{a_{j}(b)(\tau - b)^j}{j!}
				\label{eq:solution_ode_polynomial}
			\end{align}
			where each approximate derivative $a_{j}(b)$ is computed by differentiating interpolating polynomials constructed from the values in an \ODEd{} $D(r, t_n)$. For details regarding the most general formulations for the approximate derivatives $a_j(b)$ we refer the reader to \cite{buvoli2018polynomial, buvoli2019constructing}.
			
			 {\em Adams \ODEps{}} are one special family of \ODEps{} that are related to this work. Every Adams \ODEp{} can be written as
				\begin{align}
					p(\tau;b) = L_y(b) + \int_{b}^\tau L_F(\xi) d\xi
				\label{eq:adams_poly_integral_form}
				\end{align}
				where $L_y(\tau)$ and $L_F(\tau)$ are Lagrange interpolating polynomials that respectively approximate $y(t(\tau))$ and its local derivative $ry'(t(\tau))$. These polynomials can be used to construct classical Adams-Bashforth methods and their generalized block counterparts from \cite{buvoli2019constructing}.

\subsection{Parameters and notation for polynomial methods} During the timestep from $t_n$ to $t_{n+1} = t_n+h$, a polynomial method accepts inputs $y_j^{[n]}$ and produces outputs $y_j^{[n+1]}$ where ${j = 1, \ldots, q}$. Every input and output of a polynomial method approximates the differential equation solution $y(t)$ at a particular time. In local coordinates we represent these time-points using the node set $\{z_j\}_{j=1}^q$ such that
	\begin{align}
		y_j^{[n]} \approx y\left(t_n + r z_j\right) \quad \text{and} \quad y_j^{[n+1]} \approx y\left(t_n + r z_j + h\right).
	\end{align}
The input nodes of a polynomial method scale with respect to the {\em node radius} $r$, which is independent of the stepsize $h = r\alpha$. The parameter $\alpha = h/r$ is known as the {\em extrapolation factor} and its value represents the number of times the node radius divides the timestep. For large $\alpha$ the distances between a method's input nodes will be small relative to the stepsize, while the opposite is true if $\alpha$ is small. In general, polynomial methods are described in terms of $r$ and $\alpha$ rather than $r$ and $h$ since these are natural variables for working with polynomials in local coordinates.

The polynomial framework encapsulates the subset of all general linear methods \cite{butcher2006general} that are constructed using \ODEps{}. The most general formulation for a polynomial method with $s$ stages and $q$ outputs is
\begin{align}
	\begin{aligned}
		Y_i &= p_j(c_j(\alpha), b_j(\alpha)) & j &= 1, \ldots, s,\\
		y^{[n+1]}_j &= p_{j+s}(z_j + \alpha; \hspace{0.125em} b_{j+s}(\alpha)) & j &= 1, \ldots, q,
	\end{aligned}
	\label{eq:polynomial_glm}
\end{align}
where $Y_i$ denote stage values, $c_j(\alpha)$ are stage nodes, and $p_j(\tau; b)$ are all \ODEps{} constructed from an \ODEd{} containing the methods inputs, outputs, and stage values. From this general formulation it is possible to derive many different families of polynomial integrators; one example being polynomial block methods (PBMs) from \cite{buvoli2019constructing}. In this work we will generalize (\ref{eq:polynomial_glm}) for exponential integrators, and then explore the special sub-case of exponential PBMs.

\subsubsection{Propagators and iterators}
\label{subsubsec:propagator_iterator}

All polynomial integrators are parametrized in terms of the extrapolation factor $\alpha$ that scales the stepsize relative to the node radius $r$. A direct consequence of having two independent variables that control stepsize, is that all polynomial integrators naturally fall into one of two classes. If $\alpha > 0$, then the method is a {\em propagator} that advances the solution forwards in time. If $\alpha$ is zero, then the stepsize $h = r\alpha = 0$, and the method reduces to an {\em iterator} that recomputes the solution at the current time-step. %

\begin{figure}[h]
	\centering
	\subfloat[An iterator recomputes the solution at the same time points, so that $t_j^{[n+1]} = t_j^{[n]}$]{\hspace{2em}\includegraphics[height=0.09\linewidth]{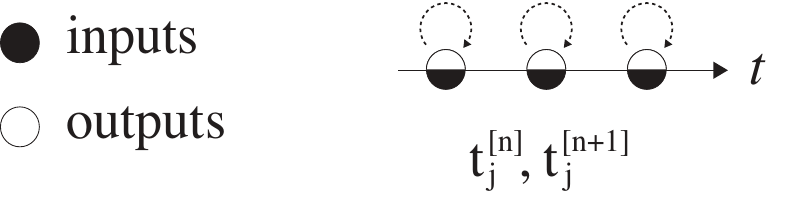} \hspace{2em}}
	\hfill
	\subfloat[A propagator advances the solution forwards in time so that  $t_j^{[n+1]} = t_j^{[n]} + r\alpha.$]{\hspace{2em}\includegraphics[height=0.09\linewidth]{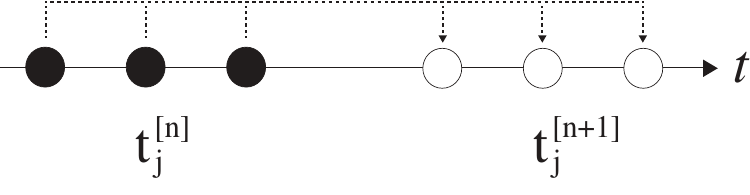}\hspace{2em}}\\
	
    \caption{Iterator and propagator methods. Input times $t_j^{[n]}$ and output times $t_j^{[n+1]}$ are shown on the real $t$ line for a method where the nodes $z_j$ are three equispaced points. Iterators can be useful for correcting or updating an approximate solution, while propagators are traditional time-stepping schemes.} 
    \label{fig:family_it_prop}
\end{figure}

In Figure \ref{fig:family_it_prop} we provide a visualization for these two types of methods. Iterators can serve several purposes, such as computing initial conditions, computing continuous output, or improving the accuracy of a candidate solution. The notion of using an iterator to improve the accuracy of a solution shares many commonalities with predictor corrector approaches for block methods \cite{shampine1969block} and spectral deferred correction methods \cite{Dutt2000SDC, christlieb2010parallel}. A detailed discussion of propagators and iterators is outside the scope of this paper. In this work, we will provide a first look at how iterators and propagators can be combined to create composite methods.

		\section{Exponential integrators}	
		\label{sec:exponential_introduction}

	Exponential integrators are a class of numerical methods for solving stiff initial value problems \cite{hochbruck2010exponentialreview}. They are derived from the prototypical semilinear equation
	 	\begin{align}
			\mathbf{y}' = \mathbf{Ly} + N(t,\mathbf{y})	&& \mathbf{y}(t_0) = \mathbf{y}_0,
			\label{eq:prototypical_semilinear}
		\end{align}
		where the autonomous linear operator $\mathbf{L}$ and the nonlinear operator $N(t,\mathbf{y})$ are chosen so that they approximate the original ODE problem. In practice, the linear operator $\mathbf{L}$ is selected so that it captures the majority of the stiffness of the right-hand side.
		
		There are two well-known classes of exponential integrators, namely {\em partitioned} exponential integrators \cite{cox2002ETDRK4, beylkin1998ELP, rainwater2014new} and {\em unpartitioned} integrators \cite{rainwater2016new, tokman2006efficient, tokman2011new}. In the following two subsections, we briefly introduce each one.
				
	\subsection{Partitioned exponential integrators}
	
	If an initial value problem can be naturally partitioned into a semilinear equation of the form (\ref{eq:prototypical_semilinear}), then it can be solved using a partitioned exponential integrator. The exact solution to (\ref{eq:prototypical_semilinear}) is obtained by applying the discrete variation of constants formula to obtain the integral equation  
	\begin{align}
		\mathbf{y}(t) = e^{(t-t_n)\mathbf{L}}\mathbf{y}(t_n) + \int^{t}_{t_n}e^{(t-s)\mathbf{L}}N(s,\mathbf{y}(s))ds.
		\label{eq:exp_derivation_formula}
	\end{align}

	Partitioned exponential integrators treat the linear term $\mathbf{L}$ exactly while approximating the nonlinearity $N(s,y(s))$ with a polynomial. For example, the second-order exponential Runge-Kutta ETDRK2 from \cite{cox2002ETDRK4}
		\begin{align}
			\begin{aligned}
			Y_1 &= e^{h\mathbf{L}} y_n + \int^{1}_0 e^{(1-s)h\mathbf{L}} N(t_n, y_n) ds, \\
			y_{n+1} &= e^{h\mathbf{L}} y_{n} + \int^{1}_0 e^{(1-s)h\mathbf{L}} \left[ N(t_n, y_n) + s(N(t_{n+1}, Y_1) - N(t_n, y_n))\right] ds,
			\end{aligned} \label{eq:erk2}	
		\end{align}
		utilizes polynomials that are constructed using stage values. Similarly, the 2nd-order exponential Adams-Bashforth method \cite{beylkin1998ELP,cox2002ETDRK4}  
		\begin{align}
			y_{n+1} &= e^{h \mathbf{L}}y_n + h\int^{1}_{0} e^{(1-s) h\mathbf{L}} \left[ N(t_n, y_{n}) + s(N(t_n, y_{n}) - N(t_{n-1}, y_{n-1}))\right]ds,	
			\label{eq:eab2}
		\end{align}
		utilizes a Lagrange interpolating polynomial constructed using previous solution values. Since exponential integrators involve weighted integrals of polynomials and exponentials, they are frequently expressed in terms of the $\varphi$-functions 
	    \begin{align}
	    	\varphi_n(z) & = 
	    	\begin{cases} 
	    		e^{z} & n = 0, \\ 
	    		\displaystyle \frac{1}{(n-1)!} \int^1_0 e^{z(1-s)} s^{n-1} ds & n> 0. 
	    	\end{cases}
	    	\label{eq:phi_function_def}
	    \end{align}
		For example, using this notation the Adams-Bashforth method (\ref{eq:eab2}) is
			\begin{align*}
				y_{n+1} &= \varphi_0(h\mathbf{L})y_n + \varphi_1(h\mathbf{L}) h N(t_n, y_{n}) + \varphi_2(hL) h(N(t_n, y_{n}) - N(t_{n-1}, y_{n-1})).	
			\end{align*}

	\subsection{Unpartitioned exponential integrators}
	\label{subsec:unpartitioned_exponential_integrators}
	
	Unpartitioned exponential integrators \cite{rainwater2016new, tokman2006efficient, tokman2011new} can be used to solve the more general initial value problem
		\begin{align}
			\mathbf{y}' = F(t, \mathbf{y}),	&& \mathbf{y}(t_0) = \mathbf{y}_0.
			\label{eq:prototypical_unpartitioned}
		\end{align}
		The key intuition is that one may obtain a localized semilinear problem of the form (\ref{eq:prototypical_semilinear}) at $t=t_n$ by rewriting the system in its autonomous form $\mathbf{y}' = F(\mathbf{y})$ and then rewriting $F(\mathbf{y})$ as

		\begin{align}
			\begin{aligned}
				F(\mathbf{y}) &= F(\mathbf{y}_n)	+ \mathbf{J}_n\left( \mathbf{y} - \mathbf{y}_n \right) + R(\mathbf{y}) \\
				R(\mathbf{y}) &= F(\mathbf{y}) - \left[ F(\mathbf{y}_n)	+ \mathbf{J}_n\left( \mathbf{y} - \mathbf{y}_n \right) \right]	
			\end{aligned}
			\label{eq:split_rhs}
		\end{align}
		where $\mathbf{y}_n = \mathbf{y}(t_n)$ and $\mathbf{J}_{n} = \frac{\partial F}{\partial y}(\mathbf{y}(t_n))$ is the Jacobian of $\mathbf{y}(t)$ at $t=t_n$. The linear operator $\mathbf{J}_n$ takes the place of $\mathbf{L}$ in (\ref{eq:prototypical_semilinear}) and the remaining terms form the nonlinearity. Given the initial condition $\mathbf{y}(t_n) = \mathbf{y}_n$, the solution of (\ref{eq:prototypical_unpartitioned}) is
			\begin{align}
				\mathbf{y}(t) = e^{(t-t_n)\mathbf{J}_n}\mathbf{y}_n + \int_{t_n}^{t} e^{(t-s)\mathbf{J}_n} \left[ F(\mathbf{y}_n)	- \mathbf{J}_n\ \mathbf{y}_n  + R(\mathbf{y})	 \right] ds.
				\label{eq:integral_unpartitioned}
			\end{align}
		 Depending on the approximation that is chosen for the remainder term $R(\mathbf{y})$ one arrives at different families of unpartitioned exponential integrators \cite{rainwater2016new, tokman2006efficient, tokman2011new}. %

\section{Exponential polynomial methods}	
\label{sec:exponential_polynomial}

	In this section we introduce the exponential equivalent of the generalized classical polynomial method (\ref{eq:polynomial_glm}), and then describe the family of exponential polynomial block methods. To accomplish this, we first extend the definitions of the \ODEd{} and the \ODEp{}. To be consistent with existing exponential integrators, the exponential extension of the \ODEp{} approximates the integral equation (\ref{eq:exp_derivation_formula}) by replacing the nonlinearity with a polynomial approximation.

	\subsection{Exponential \ODEd{} and \ODEp{}}
	
	We first discuss the extension of the \ODEd{} and \ODEp{} for the partitioned initial value problem (\ref{eq:prototypical_semilinear}). We can trivially extend the classical \ODEd{} by adding the nonlinear derivative component $N(t,\mathbf{y}(t))$ to each data element.	Since the linear term is treated exactly by an exponential integrator, there is no need to include it.

		\begin{definition}[Partitioned Exponential \ODED{}]
			An exponential \ODEd{} $\expDSP$ of size $w$ is an ordered set of tuples
			\begin{align}
				\expDSP = \left\{ \left(\tau_j,~ \mathbf{y}_j,~ r\mathbf{f}_j,~ r\mathbf{N}_j \right) \right\}_{j=1}^w	
				\label{eq:local_ode_dataset_exponential_semilinear}
			\end{align}
			where $t(\tau) = r\tau + s$, $\mathbf{y}_j \approx \mathbf{y}(t(\tau_j))$, $\mathbf{f}_j = \mathbf{L}\mathbf{y}_j + N(t(\tau_j), \mathbf{y}_j)$, and $\mathbf{N}_j = N(t(\tau_j), \mathbf{y}_j)$.
		\end{definition}

		To arrive at a generalization of the \ODEp{} (\ref{eq:solution_ode_polynomial}), we first rewrite (\ref{eq:prototypical_semilinear}) in local coordinates $\tau$, where $t(\tau) = r\tau + t_n$, and then assume that the initial condition is provided at $\tau = b$. The corresponding integral equation for the solution is
	\begin{align}
		\mathbf{y}(\tau) =  e^{(\tau-b)r\mathbf{L}} \mathbf{y}(b) 
		+ \int^{\tau}_{b} e^{(\tau - s)r\mathbf{L}} r\mathbf{N}(t(s),\mathbf{y}(t(s))
		ds.
		\label{eq:integral_equation_partitioned_local}
	\end{align}
	By expanding the nonlinearity around $\tau = b$, we can express the exact solution as
	\begin{align}
		& \mathbf{y}(\tau) =  e^{(\tau-b)r\mathbf{L}} \mathbf{y}(b) + \int^{\tau}_{b}e^{(\tau-s)r\mathbf{L}} \sum_{j=0}^\infty \frac{\mathbf{c}_{j+1}(b) (s - b)^{j}}{j!} ds, \\
		& \text{where } \quad 
		\mathbf{c}_j(b) = \left. \frac{d^{j-1}}{d\tau^{j-1}} r \mathbf{N}(t(\tau), y(t(\tau))) \right|_{\tau=b}. \nonumber
		\label{eq:integral_equation_partitioned_local_series}	
	\end{align}
	If we assume that the Taylor series converges uniformly within the domain of interest, then we can exchange the sum and the integral. Finally, using (\ref{eq:phi_function_def}) we see that the solution is an infinite series of $\varphi$-functions
	\begin{align}
		\mathbf{y}(\tau) = \varphi_0((\tau - b) r \mathbf{L}) \mathbf{y}(b) + \sum_{j=1}^\infty
				 (\tau - b)^j \varphi_j((t - b) r \mathbf{L})\mathbf{c}_j(b).
				 \label{eq:solution_phi_series}
	\end{align}
	Notice that this series is an exponential generalization of a Taylor series, since as $\mathbf{L} \to \mathbf{0}$ we recover the Taylor series for the solution expanded around $b$. 
	
	To derive the classical \ODEp{} (\ref{eq:solution_ode_polynomial}) one replaces the exact derivatives in a truncated Taylor series with polynomial approximations. Therefore, to derive its exponential equivalent, we will truncate the infinite series in (\ref{eq:solution_phi_series}) and replace the constants $\mathbf{c}_j$ with a polynomial derivative approximation of the nonlinear term. We also allow the initial condition $\mathbf{y}(b)$ to be replaced with a polynomial approximation of the solution at $\tau=b$. We formally define this type of approximant as follows.
			
	\begin{definition}[Partitioned \ppe{}]
	\label{def:polynomial_phi_expansion}
		A partitioned \ppe{} $\psi(\tau; b)$ of degree $g$ is a linear combination of $\varphi$-functions
		\begin{align}
			\psi(\tau; b) &= \varphi_0((\tau - b) r \mathbf{L}) \mathbf{a}_0(b) + \sum_{j=1}^g
				 (\tau - b)^j \varphi_j((t - b) r \mathbf{L})\mathbf{a}_j(b),			
			\label{eq:phi_expansion_phi}	
		\end{align}
		where each approximate derivative $\mathbf{a}_j$(b) must be computed using the values from a partitioned exponential ODE dataset $\expDSP = \{(\tau_j , \mathbf{y}_j , r\mathbf{f}_j, r\mathbf{N}_j)\}^w_j$ as follows:
		\begin{enumerate}
			\item The zeroth approximate derivative $\mathbf{a}_0$ is calculated in the same way as for a classical \ODEp{}. In other words, 
			\begin{align} 
				\mathbf{a}_0 = h(b),
				\label{eq:approximate_derivaties_ic}
			\end{align}
			where $h(\tau) \approx y(t(\tau))$ is a polynomial of lowest degree that interpolates at least one solution value in $\expDSP$ and whose derivative $h'(\tau)$ interpolates any number of full derivative values so that
				\begin{align*}
					&h(\tau_k) = \mathbf{y}_k && \text{for} && k \in \mathcal{A} \text{ and } \mathcal{A} \ne \emptyset,	\\
					&h'(\tau_k) = r\mathbf{f}_k && \text{for} &&  k \in B,
				\end{align*}
				and the sets $\mathcal{A}$ and $\mathcal{B}$ contain unique indices ranging from $1$ to $w$.
						
			\item The remaining approximate derivatives $\mathbf{a}_j(b)$, $j>0$, are calculated by differentiating polynomial approximations of the nonlinear derivative component $r N(\mathbf{y}(t(\tau)))$. In particular,
				\begin{align}
					\mathbf{a}_{j}(b) = \left. \frac{d^{j-1}}{d\tau^{j-1}} l_{j}(\tau) \right|_{\tau = b}, && j > 0,
					\label{eq:approximate_derivatives_nonlinearity}
				\end{align}
				where $l_j(\tau)$ is a Lagrange interpolating polynomial that interpolates at least one nonlinear derivative component value in the \ODEd{} $\expDSP$. We may express the conditions on $\l_j(\tau)$ mathematically as
					\begin{align*}
				 		\hspace{5em} l_j(\tau_k)	 &= r\mathbf{N}_k,  &&\text{for} &&  k \in \mathcal{C}^j, \text{ and } \mathcal{C}^{j} \ne \emptyset,
				 	\end{align*} 
				where $\mathcal{C}^j$ is a set that contains unique indices ranging from 1 to $w$.
		\end{enumerate}
	\end{definition}
	
	A \ppe{} generalizes a classical \ODEp{} in the sense that as $\mathbf{L} \to 0$, it reduces to (\ref{eq:solution_ode_polynomial}). To further understand the properties of \ppes{}, it is convenient to introduce a new \ODEp{} that approximates the truncated Taylor series for the local solution $rN(t(\tau), \mathbf{y}(\tau))$, expanded around the point $\tau = b$.
	
		\begin{definition}[\ODEDCP{}]
			\label{def:ode_derivative_component_polynomial}
			An exponential \ODEdcp{} of degree $g$ is a polynomial of the form
			\begin{align}
				p_N(\tau; b) &= \sum_{j=0}^{g}  \frac{\mathbf{a}_{j+1}(b)(\tau - b)^j}{j!}
				\label{eq:derivative_ode_polynomial_exponential_semilinear}
			\end{align}
			where the approximate derivatives $\mathbf{a}_j(b)$ are given by (\ref{eq:approximate_derivatives_nonlinearity}). 
 	\end{definition}
	Using this definition we can rewrite a \ppe{} in two additional ways:	
		\begin{enumerate}[leftmargin=*]
			\item The integral formulation
			\begin{align}
				\psi(\tau; b) =  e^{r\mathbf{L}(\tau-\tau_0)}
				\underbracket[.03em][0.2em]
				{
					\mathbf{a}_0(b)
				}_{
					\begin{subarray}{c} 
						| \\
						\mathclap{\text{Interpolating polynomial approximating $y(t(\tau))$ evaluated at $b$.}}
					\end{subarray}
				}
				+ \int^{\tau}_{b} e^{r\mathbf{L}(\tau - s)} 
				\overbracket[.03em][0.3em]
				{
					p_N(s; b)
				}^{
					\begin{subarray}{c} 
						\mathclap{\text{\ODEdcp{} expanded at $b$.}} \\ 
						|
					\end{subarray}
				}
				ds.
				\label{eq:phi_expansion_integral}
			\end{align}
			This representation explicitly shows the polynomial approximations that have been substituted into the integral equation (\ref{eq:integral_equation_partitioned_local}). To show equivalence with the $\varphi$-expansion (\ref{eq:phi_expansion_phi}), apply the change of variables $s = (\tau - b) \nu + b$, and simplify using (\ref{eq:phi_function_def}).	\\[-0.5em]	

			\item The differential formulation; $\psi(\tau; b)$ is the solution of the linear system 
			\begin{align}
				\begin{aligned}
					\mathbf{y}'(\tau) &= r\mathbf{L}\mathbf{y} + p_N(\tau; b), && \mathbf{y}(b) = \mathbf{a}_0(b).				
				\end{aligned}
				\label{eq:phi_expansion_differential}
			\end{align}
			From this representation we see that all \ppes{} can be computed by solving a linear differential equation, and that it is possible to use the value of a \ppe{} at time $t_1$ as an initial condition for computing the \ppe{} at time $t_2 > t_1$. To show equivalence to ($\ref{eq:phi_expansion_integral}$) simply apply the integrating factor method with $e^{\mathbf{L}\tau}$. \\[-0.5em]
		\end{enumerate}

	As is the case for a classical \ODEp{}, there are many ways to select the approximate derivatives for a \ppe{}. To simplify method construction, we introduce a special family of polynomial $\varphi$-expansions that is related to the Adams \ODEp{} (\ref{eq:adams_poly_integral_form}) for classical integrators.
 
	\begin{definition}[\APE{}]
	\label{def:adams_phi_expansion}	
	 The polynomial $\varphi$-expansion (\ref{eq:phi_expansion_integral}, \ref{eq:phi_expansion_differential}, \ref{eq:phi_expansion_phi}) is of Adams type if it satisfies the following two conditions:
		\begin{enumerate}
	 		\item The interpolating polynomial $h(\tau)$, used to set the zeroth approximate derivative $\mathbf{a}_0(b) = h(b)$, is a Lagrange interpolating polynomial $L_y(\tau)$ that interpolates at least one solution value in an \ODEd{}.
			\item The \ODEdcp{} $p_N(\tau; b) = L_N(\tau)$ where $L_N(\tau)$ is a Lagrange interpolating polynomial that interpolates any number of derivative component values in an \ODEd{}.	
	 	\end{enumerate}
	 \end{definition}
		The three equivalent representations for an Adams \ppe{} are:
		\begin{small}
		\begin{align*}
			&\text{1. integral} & \psi(\tau; b) &= e^{r\mathbf{L}(\tau-b)} L_y(b) + \int^{\tau}_{b} e^{r\mathbf{L}(\tau - s)}L_N(s) 	ds,\\[0.5em]
			&\text{2. differential} & 	\mathbf{y}'(\tau) &= r\mathbf{L}\mathbf{y} + L_N(\tau), \quad \mathbf{y}(b) = L_y(b),	\\[0.5em]
			&\text{3. $\varphi$-expansion} & \psi(\tau; b) &= \varphi_0(r\mathbf{L}(\tau - b))L_y(b) + \sum_{k=0}^g (\tau - b)^{k+1} L_N^{(k)}(b) \varphi_{k+1}(r\mathbf{L}(\tau - b)). 
		\end{align*}
		\end{small}%
	We name these simplified $\varphi$-expansions \Apes{} for two reasons. First, they reduce to a classical Adams \ODEp{} (\ref{eq:adams_poly_integral_form}) as $\mathbf{L}\to \mathbf{0}$. Second, we can derive exponential Adams-Bashforth and Adams-Moulton methods \cite{beylkin1998ELP} using an \Ape{} constructed from an \ODEd{} with equispaced nodes.

	\subsubsection{Obtaining coefficients for \ppes{}}
	
	From (\ref{eq:phi_expansion_phi}) we see that all \ppes{} are linear combinations of $\varphi$-functions where the weights are expressed in terms of the approximate derivatives $\mathbf{a}_j(b)$. All approximate derivatives can be written as linear combinations of the data elements in the corresponding exponential \ODEd{}, and a detailed procedure for obtaining these weights is described in \cite[Sec 3.7]{buvoli2018polynomial}. %
		For \Apes{}, the procedure for obtaining coefficients is simpler and is described in Appendix \ref{ap:adams_ppe_coefficients}.
		
		\subsubsection{Extension to unpartitioned problems} We can also extend the \ODEd{} and the \ODEp{}  for the unpartitioned problem (\ref{eq:prototypical_unpartitioned}) using local linearization. 
	To derive an unpartitioned \ppe{} we
		\begin{enumerate}
			\item rewrite (\ref{eq:prototypical_unpartitioned}) in autonomous form,
			\item re-express the equation in local coordinates $\tau$,
			\item locally linearize around $\tau = \ell$ to obtain the equation (\ref{eq:split_rhs}) with ${\mathbf{y}_n \to \mathbf{y}(\ell)}$ and $\mathbf{J}_n \to \frac{\partial F}{\partial y}(\mathbf{y}_\ell)$, and
			
			\item assume that the initial condition is supplied at $\tau = b$. 
		\end{enumerate}
	Using the integrating factor method, the corresponding solution is 
			\begin{align}
				\mathbf{y}(\tau) = e^{(\tau-b)r\mathbf{J}(\mathbf{y}_\ell)}\mathbf{y}_0 + r \int_{b}^{\tau} e^{(\tau-s)r\mathbf{J}(\mathbf{y}_\ell)} \left[ F(\mathbf{y}_\ell)	- \mathbf{J}(\mathbf{y}_\ell)\ \mathbf{y}_\ell  + R(\mathbf{y})	 \right] ds,
				\label{eq:upartitioned_solution}
			\end{align}
			where $\mathbf{y}_0 = \mathbf{y}(\tau_0)$, $\mathbf{y}_\ell = \mathbf{y}(\ell)$, and $\mathbf{J}(y) = \frac{\partial F}{\partial y}(y)$ is the Jacobian of the right-hand-side evaluated at $y$. Next we introduce the {\em unpartitioned \ODEd{}}
			\begin{align}
				\expDSU = \left\{ \left(\tau_j,~ \mathbf{y}_j,~ r\mathbf{f}_j, r\mathbf{R}_j \right) \right\}_{j=1}^w, && \mathbf{R}_j = R(\mathbf{y}_j),
				\label{eq:upartitioned_exponential_ode_dataset}		
			\end{align}
			and an associated \ODEdcp{} $p_R(\tau; b)$ that is defined identically to the \ODEdcp{} in Definition \ref{def:ode_derivative_component_polynomial}, except that the approximate derivatives are now formed using $\mathbf{R}_j$ instead of $\mathbf{N}_j$. To obtain an unpartitioned \ppe{}, we make the following replacements to (\ref{eq:upartitioned_solution}):
				\begin{align}
					\mathbf{y}_0 \to  \mathbf{a}_0(b), \quad
					\mathbf{y}_l \to  \hat{\mathbf{a}}_0(\ell), \quad
					\text{and} \quad
					R(\mathbf{y}(s))  \to p_R(s; b),
					\label{eq:integral_unpartitioned_substitutions}
				\end{align}
				where $\mathbf{a}_0$ and $\hat{\mathbf{a}}_0$ are zeroth-order approximate derivatives defined in (\ref{eq:approximate_derivaties_ic}), and $p_R(s; b)$ is an unpartitioned \ODEdp{} expanded at $b$. 
				
				In comparison to the partitioned \ppe{} we have an additional parameter $\ell$ which specified the temporal location where the unpartitioned ODE was locally linearized. Once again, we can limit the number of free parameters by introducing a special family of Adams approximations.

		\begin{definition}[Unpartitioned \APE{}] Let $L_y(\tau) \approx \mathbf{y}(t(\tau))$ and $L_R(\tau) \approx rR(\mathbf{y}(t(\tau)))$ be two Lagrange interpolating polynomials constructed from the \ODEd{} (\ref{eq:upartitioned_exponential_ode_dataset}) such that $L_y(\tau)$ interpolates at least one solution value, and $L_R(\tau)$ interpolates at any number of derivative component values;
		\begin{align}
				L_y(\tau_k) = \mathbf{y}_k, \text{ for } k \in \mathcal{A} \text{ where } |\mathcal{A}| \ge 1,  \text{ and } L_R(\tau_j) = r R(\mathbf{y}_k), \text{ for } k \in \mathcal{B}.
		\end{align}
		An unpartitioned \Ape{} is an approximation of the form (\ref{eq:upartitioned_solution}, \ref{eq:integral_unpartitioned_substitutions}) where the expansion point $\ell$ and left integration point $b$ are equal such that $\ell = b$, and $\mathbf{a}_0(b) = \hat{\mathbf{a}}_0(b) = L_y(b)$ and $p_R(s; b) = L_R(s)$.
		\end{definition}
		If we let $\mathbf{J}_b = \frac{\partial F}{\partial y}(L_y(b))$ then we can write the formula for an unpartitioned Adams \ppe{} in the usual three ways: 
		\begin{small}
		\begin{align*}
			&\text{differential} & \mathbf{y}'(\tau) &= rF(L_y(b)) + r\mathbf{J}_b(\mathbf{y}(\tau) - L_y(b)) + L_R(\tau), \quad \mathbf{y}(b) = L_y(b),\\[0.5em]
			&\text{integral} & \hat{\psi}(\tau; b) &= e^{r\mathbf{L}(\tau - b)} L_y(b) + \int_{b}^{\tau} e^{(\tau-s)r\mathbf{J}_b} \left[ rF(L_y(b))	- r\mathbf{J}_b L_y(b)  + L_R(s)	 \right] ds,\\[0.5em]
			& & &= L_y(b) + \int_{b}^{\tau} e^{(\tau-s)r\mathbf{J}_b} \left[ rF(L_y(b))  + L_R(s)	 \right] ds,\\[0.5em]
			&\text{$\varphi$-expansion} & \hat{\psi}(\tau; b) &= 
				\begin{aligned}[t]
					& L_y(b) + (\tau - b)\varphi_1(r\mathbf{J}_b(\tau - b)) rF(L_y(b))  \\
					& + \sum_{k=0}^g (\tau - b)^{k+1} L^{(k)}_R(b) \varphi_{k+1}(r\mathbf{J}_b(\tau - b)).
				\end{aligned}			
		\end{align*}
		\end{small}%
		The identity $e^{h\mathbf{A}}\mathbf{x} = \mathbf{x} + h\varphi_1(h \mathbf{A})\mathbf{A}\mathbf{x}$ was used to write the $\varphi$-expansion and the second formula for the integral formulation.

	\subsection{A general formulation for exponential polynomial methods}
	
	The general formulation for a partitioned exponential polynomial integrator with $s$ stages and $q$ outputs is
	\begin{align}
		\begin{aligned}
			Y_i &= \psi_j(c_j(\alpha);~ b_j(\alpha)) & j &= 1, \ldots, s,\\
			y^{[n+1]}_j &= \psi_{j+s}(z_j + \alpha; \hspace{0.125em} b_{j+s}(\alpha)) & j &= 1, \ldots, q,
		\end{aligned}
		\label{eq:exponential_glm}
	\end{align}
	where $Y_i$ are stage values and $\psi_j(\tau; b)$ are partitioned \ppes{} that are constructed from a partitioned exponential \ODEd{} that contains the methods inputs, outputs, and stage values. The definition for an unpartitioned exponential polynomial integrator is identical except we replace $\psi_j$, $j=1, \ldots, q+s$, with unpartitioned \ppes{} $\hat{\psi}_j$ that also depend on the parameters $l_j(\alpha)$.
	
	Eqn. (\ref{eq:exponential_glm}) encapsulates all families of polynomial exponential integrators, and is the exponential generalization of (\ref{eq:polynomial_glm}). However, the space of such methods is too large to explore in this initial work. We therefore restrict ourselves to the family of exponential polynomial block methods.
	
	\subsection{Exponential polynomial block methods}
	
	Classical block methods \cite{shampine1969block} are multivalued integrators that advance multiple solution values per timestep. Depending on the structure of their coefficient matrices they compute their outputs either in serial or in parallel \cite{gear1988parallel, gander201550}. Block methods provide a good starting point for deriving high-order polynomial integrators, since they have no stage values and their multiple inputs can be used to easily construct high-order polynomial approximations. Parallel polynomial block methods with nodes in the complex plane were introduced in \cite{buvoli2019constructing}.
 
	As is the case for classical PBMs, exponential polynomial block methods (EPBMs) are a good starting point for deriving high-order exponential polynomial integrators. A partitioned EPBM depends on the parameters:
	\begin{center}
			\renewcommand*{\arraystretch}{1.5}
			\begin{tabular}{llllll}
				$q$ 		& number of inputs/outputs 	& \hspace{3em} 	& $\left\{ z_j \right\}_{j=1}^q$ & nodes, $z_j \in \mathbb{C}$, $|z_j| \le 1$ \\
				$r$ 		& node radius, $r \ge 0$ 	& 				& $\left\{ b_j \right\}_{j=1}^q$ & expansion points \\
				$\alpha$ 	& extrapolation factor 		& 				& \\
			\end{tabular}
	\end{center}
	and can be written as	
	\begin{align}
		y^{[n+1]}_j &= \psi_j(z_j + \alpha;~ b_j(\alpha)), & j & =1, \ldots,q,
		\label{eq:polynomial_block_method_partitioned}
	\end{align}
	where each $\psi_j(\tau; b)$ is an \ppe{} built from the exponential \ODEd{} 
	\begin{align}
		\expDSP[t_n] = \underbrace{\left\{ \left(z_j,~ \mathbf{y}_j^{[n]},~ r\mathbf{N}_j^{[n]} \right) \right\}_{j=1}^q}_{\text{method inputs}} \bigcup \underbrace{\left\{ \left(z_j+\alpha,~  \mathbf{y}_j^{[n+1]},~ r\mathbf{N}_j^{[n+1]}\right) \right\}_{j=1}^q}_{\text{method outputs}}.
	\end{align}	
	We can also write any partitioned EPBM in coefficient form as
	\begin{align}
		\begin{aligned}
			y^{[n+1]}_j &= \varphi_0(\eta_j \mathbf{L}) \sum_{k=1}^q \left(A_{jk} \mathbf{y}_k^{[n]} + C_{jk} \mathbf{y}_k^{[n+1]} \right) \\
			&~ + r \sum_{k=1}^{q} \varphi_k(\eta_j \mathbf{L}) \sum_{l=1}^q \left( B_{jkl} \mathbf{N}^{[n]}_l + D_{jkl} \mathbf{N}^{[n+1]}_{l}\right),
		\end{aligned}
	\end{align}
	where $\eta_j = (z_j + \alpha - b_j)$ and $j = 1, \ldots, q$. Similarly we can write an unpartitioned EPBM as 
	\begin{align}
		y_j^{[n+1]} = \hat{\psi}_j(z_n + \alpha;~ b_j(\alpha)), \quad j = 1, \ldots, q,	
		\label{eq:polynomial_block_method_unpartitioned}
	\end{align}
	or in its coefficient form
	\begin{align}
		\begin{aligned}	
				y^{[n+1]}_j &= \left(A_{jk} \mathbf{y}_k^{[n]} + C_{jk} \mathbf{y}_k^{[n+1]} \right) + \varphi_1(\eta_j \mathbf{J}_{b_j}) F\left(\sum_{k=1}^q \left(A_{jk} \mathbf{y}_k^{[n]} + C_{jk} \mathbf{y}_k^{[n+1]} \right)\right) \\
			&~ + r \sum_{k=1}^{q} \varphi_k(\eta_j \mathbf{J}_{b_j}) \sum_{l=1}^q \left( B_{jkl} \mathbf{R}^{[n]}_l + D_{jkl} \mathbf{R}^{[n+1]}_{l}\right).
		\end{aligned}	
	\end{align}

		\section{Constructing Adams exponential polynomial methods}
		\label{sec:constructing_pbm}
		
		In this section we discuss several approaches for constructing EPBMs with Adams \ppes{}. In particular, we derive several new classes of both parallel and serial EPBMs, and also discuss a strategy for obtaining initial conditions by composing multiple iterator methods. %

	\subsection{Parameter selection}
	\label{subsec:constructing_epbm}
	To simplify the construction of Adams EPBMs it is convenient to rewrite them in the integral form
	\begin{align}
		&\text{Partitioned:} \nonumber \\
		& y^{[n+1]}_j = e^{r\eta_j\mathbf{L}} L^{[j]}_y(b_j) + \int^{z_j + \alpha}_{b_j} e^{r\mathbf{L}(\tau - s)}L^{[j]}_N(s) ds, \label{eq:epbm_partitioned_output} \\ 
		&\text{Unpartitioned:} \nonumber \\
		& y^{[n+1]}_j = L^{[j]}_y(b_j) + r\eta_j \varphi_1(r\eta_j \mathbf{J}_{b_j}) F(L^{[j]}_y(b_j)) + \int^{z_j + \alpha}_{b_j} e^{r\mathbf{J}_{b_j}(\tau - s)}L^{[j]}_R(s) ds, \label{eq:epbm_unpartitioned_output}
	\end{align}
	where $\eta_j = z_j + \alpha - b_j$ and $j = 1, \ldots, q$. To construct these types of methods we must select:
	\begin{enumerate}
		\item a set of nodes $\left\{ z_j \right\}_{j=1}^q$ that determine the input and output times,
		\item the Lagrange polynomials $L_y^{[j]}$, $L_N^{[j]}$ or $L_y^{[j]}$, $L_R^{[j]}$ that respectively interpolate values or derivatives at the input and output nodes, and
		\item the integration endpoints $b_j$.
	\end{enumerate}
	 In the following subsections we describe multiple strategies for selecting each of these parameters. However, we note that the parameters for exponential methods can also be selected in entirely different ways and in a different order than what is shown in this paper.
	
	\subsubsection{Node selection}
		
		Polynomial methods can be constructed using either real-valued or complex-valued nodes $\left\{ z_j \right\}_{j=1}^q$. In general, the node type has a nontrivial effect on the linear stability of a method. For example, for diagonally implicit polynomial methods, imaginary equispaced nodes offer improved stability compared to real-valued equispaced nodes \cite{buvoli2019constructing}. For serial methods, the node ordering will also affect stability since it determines the output ordering.
		
		With the exception of Section \ref{subsubsection:imaginary-nodes}, we only consider real-valued nodes that are normalized so they lie on the interval $[-1,1]$, and are ordered so that $$z_1 < z_2 < \ldots < z_q.$$
		This is of course only one possible choice, however a complete discussion on optimal node selection is nontrivial, and lies outside the scope of this paper.		
		
	\subsubsection{Selecting the polynomial $L_y^{[j]}(\tau)$ and the endpoints $b_j$}
	
	The polynomial $L_y^{[j]}(\tau)$ approximates the initial condition for the integral equation (\ref{eq:integral_equation_partitioned_local}) and  (\ref{eq:upartitioned_solution}). To construct high-order parallel methods, we choose the highest-order polynomial that does not lead to a method with coupled outputs.  This property is obtained if we let ${L_y^{[j]}(\tau) = L_y(\tau)}$, where $L_y(\tau)$ is a Lagrange interpolating polynomial of order $q-1$ that interpolates through each of the methods inputs so that
		\begin{align}
			L_y(\tau) = \sum_{j=1}^q \ell_j(\tau) y_j^{[n]}, \quad \text{where}\quad \ell_j(\tau)=\prod_{\substack{k = 1\\k \ne j }}^q \frac{\tau - z_k}{z_j - z_k}.
			\label{eq:lagrange_Ly_formula}
		\end{align}
		To simplify things even further, we select the integration endpoints so that they are equal to any of the node values; in other words, there must exist an integer $k(j)$ such that $b_j = z_{k(j)}$ for $j = 1, \ldots, q$. Under these assumptions, the formulas (\ref{eq:epbm_partitioned_output}) and (\ref{eq:epbm_unpartitioned_output}) for the $j$th output of an EPBM reduce to
		\begin{align}
			&\text{partitioned:} \nonumber \\
			& y^{[n+1]}_j = e^{r\eta_j \mathbf{L}} y^{[n]}_{k(j)} + \int^{z_j + \alpha}_{b_j} e^{r\mathbf{L}(\tau - s)}L^{[j]}_N(s) ds, \label{eq:adams_pbm_output_partitioned} \\
			&\text{unpartitioned:} \nonumber \\
			& y^{[n+1]}_j = y^{[n]}_{k(j)} + r\eta_j \varphi_1(r\eta_j \mathbf{J}_{b_j}) F(y^{[n]}_{k(j)}) + \int^{z_j + \alpha}_{b_j} e^{r\mathbf{J}_{b_j}(\tau - s)}L^{[j]}_R(s) ds.	 \label{eq:adams_pbm_output_unpartitioned}
		\end{align}

	\subsubsection{Selecting the polynomials $L_N^{[j]}(\tau)$ and $L_R^{[j]}(\tau)$}
	\label{subsubsub:Lagrange_N}
	
\begin{figure}[h]
			\centering
			\includegraphics[width=1\linewidth]{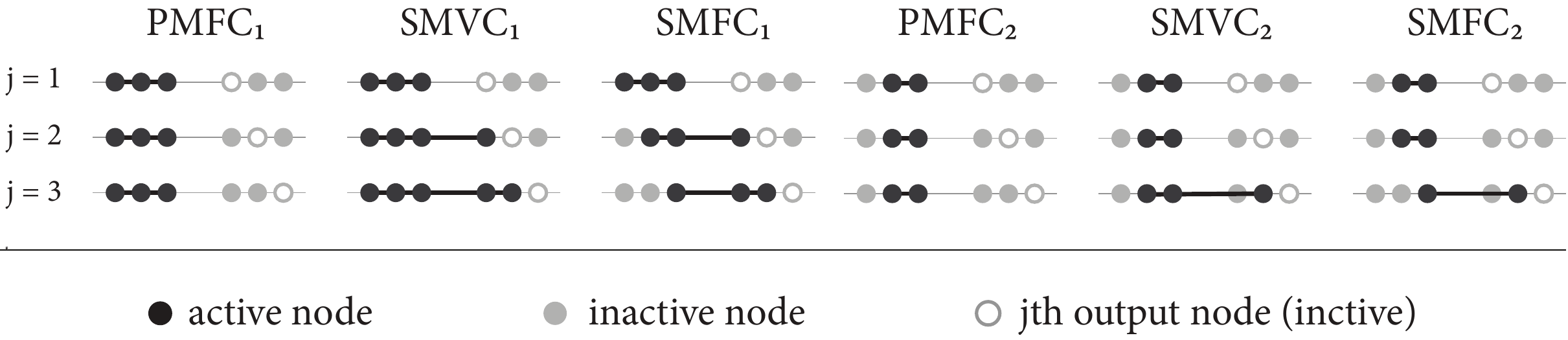}	
			\caption{A stencil-like diagram showcasing the three construction strategies for methods with three real-valued outputs ($q=3$) and equispaced nodes $\{z_j\}_{j=1}^3=\{-1,0,1\}$. Each type of method must compute three outputs and therefore requires three polynomials. The stencils for each of these polynomials are stacked vertically. The horizontal gray lines represent the time axis that flows to the right. The leftmost group of three dots on each time-line represents the input nodes while the rightmost group of three dots represent the output nodes. The diagrams show the temporal locations of the data that is used to construct the Lagrange polynomials that approximate the nonlinear term. A node is colored black and labeled {\em active} if the polynomial $L^{[j]}(\tau)$ interpolates data at this node. Conversely, the node is colored light gray and labeled {\em inactive} if the data at the node is unused. The inactive output node $z_j + \alpha$, which is the location where the \ppe{} is evaluated, is shown as a white circle with a gray border.}
			\label{fig:active_node_diagrams_real_explicit}
			
		\end{figure}

	The Lagrange polynomials $L_N^{[j]}(\tau)$ and $L_R^{[j]}(\tau)$ respectively approximate the nonlinear terms in the integral equations (\ref{eq:integral_equation_partitioned_local}) and (\ref{eq:upartitioned_solution}). Here we propose three strategies for choosing these polynomials that lead to either parallel or serial methods. %
 We will first describe the motivation behind our choices and then present formulas that determine the indices of the data used to form the Lagrange polynomials. In Figure \ref{fig:active_node_diagrams_real_explicit} we also present several example stencils that graphically illustrate the temporal locations of the data. To appreciate the simple geometric construction behind each of the proposed polynomials, the descriptions and formulas presented below should be read in tandem with their corresponding diagrams in Figure \ref{fig:active_node_diagrams_real_explicit}.
 
 	Before describing our choices, we introduce the {\em input index set} $I(j)$ and the {\em output index set} $O(j)$ for the Lagrange polynomials $L_N^{[j]}(\tau)$ and $L_R^{[j]}(\tau)$. These sets respectively contain all of the indices of the input data and the output data that is used to construct the polynomial. For example, if $L_N^{[j]}(\tau)$ or $L_R^{[j]}(\tau)$ has an interpolation node at an input node $\tau = \tau_k$, then $k$ would be a member of $I(j)$. We formally define the sets as:
		 \begin{center}
		 	\vspace{0.25em}
		 	\renewcommand{\arraystretch}{1.5}
		 	\begin{tabular}{ll}
		 		partitioned & unpartitioned \\ 
		 		$k \in I(j) \hspace{0.25em} \iff  L^{[j]}_N(z_k) = N_k^{[n]}$, \hspace{1.5em}  & $k \in I(j) \hspace{0.25em} \iff L^{[j]}_R(z_k) = R_k^{[n]}$,	 \\
		 		$k \in O(j) \iff  L^{[j]}_N(z_k+\alpha) = N_k^{[n+1]}$, \hspace{1.5em}  & $k \in O(j) \iff L^{[j]}_R(z_k+\alpha) = R_k^{[n+1]}$.	
		 	\end{tabular}
		 	\vspace{.25em}
		 \end{center}
		The order of the polynomial $L^{[j]}_N(z_k)$ or $L^{[j]}_R(z_k)$ is therefore $\left|I(j)\right| + \left|O(j)\right| - 1$. 
		
		Our three proposed strategies for the construction $L_N^{[j]}(\tau)$ and $L_R^{[j]}(\tau)$ are:
		\begin{enumerate}[leftmargin=*]
			\item {\em Parallel maximal-fixed-cardinality$-\ell$} ({\bf PMFC$_\ell$}): The Lagrange polynomial interpolates the nonlinear derivatives at all input nodes with index greater than or equal to $\ell$. The set $O(j)$ must be empty to retain parallelism.
			
			\item {\em Serial maximal-variable-cardinality-$\ell$} ({\bf SMVC$_\ell$}): The Lagrange polynomial interpolates the nonlinear derivatives at all input nodes and previously computed outputs with index greater than or equal to $\ell$. The cardinality of set $O(j)$ grows as we add new information. 			
			
			\item {\em Serial maximal-fixed-cardinality-$\ell$} ({\bf SMFC$_\ell$}): The Lagrange polynomial interpolates the nonlinear derivatives at  all previously computed outputs and some of the inputs with index greater than or equal to $\ell$. To keep the cardinality of $I(j) \cup O(j)$ fixed for all $j$, we drop inputs in favor of more recently computed outputs. 	
		\end{enumerate}
		The formulas for the three construction strategies are written in Table \ref{tab:active_index_sets_real}.
		Finally, implicit EPBMs can be created by taking $O(j) \to O(j) \cup {j}$; however, a detailed discussion of these methods is outside the scope of this paper.
		
		\begin{table}[h]
			\centering
			\renewcommand{\arraystretch}{1.5}
			\begin{tabular}{r|ll}
			& $I(j)$ - input index set & $O(j)$  - output index set\\ \hline
			{\bf PMFC$_\ell$} & $\left\{ \ell, ~\ell+1~ \ldots,~q\right\}$ & $\left\{  \right\}$ \\
			{\bf SMFC$_\ell$} & $\left\{ \max(j,\ell),~\max(j,\ell)+1,~ \ldots,~q\right\}$ & $\left\{ \ell,~\ell+1,~ \ldots,~ j-1 \right\}$ \\
			{\bf SMVC$_\ell$} & $\left\{ \ell,~\ell+1,~ \ldots,~q\right\}$ & $\left\{ \ell,~\ell+1,~ \ldots,~ j-1 \right\}$
			\end{tabular}
			\caption{Formulas describing the input index sets $I(h)$ and the output index sets $O(j)$ for three construction strategies.}
			\label{tab:active_index_sets_real}
		\end{table}
		
		\subsection{Iterator methods for obtaining initial conditions}
				
				EPBMs require multiple inputs at the first time-step. Since we are normally provided with only one initial value, an exponential Runge-Kutta scheme can be used to compute the solution at each of the starting nodes. However, this approach requires implementing a separate method. Moreover, it may not always be possible to match the order of a starting method with that of the EPBM. Here we present an alternative strategy for obtaining initial conditions by repeatedly applying an iterator method (i.e. a polynomial method with $\alpha = 0$, described in Section \ref{subsubsec:propagator_iterator}). 
				
				If a set of input values contains at least one accurate initial condition, then it is possible to construct a polynomial method that iteratively improves the accuracy of the other solution values \cite[Chpt. 6.2]{buvoli2018polynomial}. This is accomplished by repeatedly applying a polynomial iterator method that uses the accurate value as the initial condition in a discrete Picard iteration. Each application of the method increases the order-of-accuracy of the inputs by one. The maximum achievable accuracy is bounded above by the order of the iterator or the order of the accurate initial condition. This idea is closely related to spectral deferred correction methods \cite{Dutt2000SDC}, except we do not use an error formulation and we do not discard any values after the iteration is complete.
				
				Here we describe two exponential iterators that can be used to produce input values for any EPBM. They are both constructed using the strategies from Section \ref{subsec:constructing_epbm}. For simplicity we assume that the initial condition is given at the first node so that $\mathbf{y}_1^{[0]} = y(rz_1 + t_0)$ (if this is not true, then simply redefine $t_0$). We first obtain a coarse estimate for the remaining initial values by assuming that the solution is equal at all the nodes. We then repeatedly apply an iterator to improve the accuracy of this coarse estimate. We can abstractly write the formula for this iteration as
				\begin{align}
						\mathbf{y}^{[0]} = M^k \mathbf{c}, && \mathbf{c}_j = \mathbf{y}_1^{[0]}, \quad j = 1, \ldots q,
				\end{align}
				where $\mathbf{c}$ denotes the initial coarse approximation, $M$ denotes the PBM iterator method, and the computation involves applying the method $k$ times to the course initial condition $\mathbf{c}$.
					
				 To construct the method $M$ we propose the following parameters: $L_y^{[j]}(\tau) = L_y(\tau)$ from (\ref{eq:lagrange_Ly_formula}), $b_j = z_1$, and $L_N^{[j]}(\tau)$ \text{ or } $L_R^{[j]}(\tau)$ should be constructed in accordance with the formulas for PMFC$_\ell$ or SMFC$_\ell$. Selecting PMFC$_\ell$ leads to a parallel iterator while SMFC$_\ell$ leads to a serial iterator. The SMVC$_\ell$ construction cannot be used since the input and output nodes overlap when $\alpha = 0$. For typical node choices $\{ z_j \}_{j=1}^q$ it is sufficient to run the iteration $q$ times; however, node sets that allow for higher accuracy can benefit from additional iterations. Finally, as $k\to\infty$ the iteration converges to the stages and outputs of a fully implicit exponential collocation method.
				 
				Iterators can also be applied to construct composite EPBMs. An iterator and a propagator that share the same nodes can be combined to create a composite method whose computation involves first advancing the timestep with the propagator and then correcting the output $\kappa$ times with the iterator. The composite method can be written abstractly as
				\begin{align}
					\mathbf{y}^{[n+1]} = M^\kappa P \mathbf{y}^{[n]}	,
				\end{align}
				where $P$ denotes a propagator and $M$ the iterator. As we show in our numerical experiments, this method construction can lead to methods with improved stability properties.

	\subsection{Parallel Adams PBMs with Legendre nodes}
	
	The parameter choices proposed in Section \ref{subsec:constructing_epbm} can be used to construct  families of partitioned or unpartitioned exponential methods. To obtain one single method it is necessary to select a set of nodes and the extrapolation parameter $\alpha$. 
	
	For the sake of brevity, we only analyze parallel polynomial exponential methods where $L_N(\tau)$ is constructed in accordance to the PMFC strategy. Furthermore, we only consider real-valued nodes that are given by the union of negative one and the $q-1$ Legendre points so that
		\begin{align}
			\{z_j\}_{j=1}^q = \{-1, x_1, \ldots, x_{q-1}\} 
			\label{eq:leg_zero_nodes}
		\end{align}
		where $x_j$ is the $j$th zero of the $q$th Legendre polynomial $P_{q}(x)$. Since we are using Legendre points we choose the PMFC$_2$ construction strategy so that the polynomials for $L^{[j]}_N(\tau)$ or $L^{[j]}_R(\tau)$ are both constructed using only the Legendre nodes. We also select $b_j = z_1$ so that the integration endpoint is negative one. We describe the methods parameters in Table \ref{tab:leg_methods}.
		\begin{table}[h!]
			\centering
			\renewcommand{\arraystretch}{1.5}
			\begin{tabular}{rcl}
				$\{z_j\}_{j=1}^q$ &:& $\{-1, x_{q-1}, \ldots, x_{q-1}\}$ from Eqn. (\ref{eq:leg_zero_nodes}) \\  
				$L^{[j]}_y(\tau)$ &:& $\sum_{j=1}^q \ell_j(\tau) y_j^{[n]}$ from Eqn. (\ref{eq:lagrange_Ly_formula}) \\  
				$L^{[j]}_N(\tau)$ or $L^{[j]}_R(\tau)$ &:& PMFC$_2$ \\
				$\{ b_j \}_{j=1}^q$ &:& $\{ z_1 \}_{j=1}^q$
			\end{tabular}
			\caption{Parameters for Legendre EPBMs. The coefficients for methods with $q=2,3,4$ are shown in Appendix \ref{ap:leg_method_coeff}.}	
			\label{tab:leg_methods}
		\end{table}
		
		Like all polynomial methods, these Legendre EPBMs are parametrized in terms of the extrapolation factor $\alpha$; we can write the method abstractly as $M(\alpha)$. Different choices will impact both the method's stability, accuracy, and susceptibility to round-off errors. We will primarily focus on Legendre propagators with $\alpha = 2$, and Legendre iterators with $\alpha = 0$ for computing initial conditions. For dispersive equations we also briefly consider composite methods of the form
			\begin{align}
				M_{\text{composite}}(\alpha) = M(0)^\kappa M(\alpha).
				\label{eq:leg_method_composite}
			\end{align}

		In Table \ref{tab:EPBM_pseudocode} we show pseudocode for the composite method with Legendre nodes. If $\kappa$=0, then the method reduces to the standard Legendre EPBM.
		
		Finally, though we have found that Legendre based nodes lead to efficient methods, they are by no means the optimal choice. However, optimal node selection is nontrivial and should be treated as a separate topic.

		\begin{table}[t]
\renewcommand{\arraystretch}{1.25}
\begin{tabularx}{\textwidth}{|X|}
  \hline
  
  \vspace{-0.75em}
  {\bf Partitioned Composite EPBM (\ref{eq:leg_method_composite}) with $M(\alpha)$ from Table \ref{tab:leg_methods} } \\[.25em]   \hline      
        \vspace{-0.5em}
        
        {\color{gray} \# -- Definitions -------------------------------------------------------------------------------------- } \\
        $\eta_j(\alpha) = z_j + \alpha + 1$ \\
        $ \hat{\ell}_j(\tau) = \prod_{k=1,k\ne j}^{q-1} \frac{\tau - x_k}{x_j - x_k}$ {\color{gray} \# $x_j$ are the Legendre nodes; see Eqn. (\ref{eq:leg_zero_nodes}).}\\
        $L_N(\tau) = \sum_{j=1}^{q-1} \hat{\ell}_j(\tau) \mathbf{N}_{j+1}$ \\
        
        {\color{gray} \# -- Propagator ------------------------------------------------------------------------------------} \\
        \text{{\bf for} j = 2 to $q$} {\color{gray} \# {\bf parallelizable loop}} \\        
        \hspace{1.5em} $\mathbf{N}_j = N(y_j^{[n]})$
        
        \text{{\bf for} j = 1 to $q$} {\color{gray} \# {\bf parallelizable loop}} \\        
        \hspace{1.5em} 
        	$y_j^{[n+1]} = \varphi_0(r\mathbf{L}\eta_j(\alpha))y_1^{[n]} + \sum_{k=0}^{q-1} (\eta_j(\alpha))^{k+1} \varphi_{k+1}(r\mathbf{L}\eta_j(\alpha)) L_N^{(k)}(-1) $ \\
        {\color{gray} \# -- Iterator ----------------------------------------------------------------------------------------- } \\
        \text{{\bf for} k = 1 to $\kappa$} \\ 
        \hspace{1.5em} \text{{\bf for} j = 2 to $q$} {\color{gray} \# {\bf parallelizable loop}} \\        
        \hspace{3em} $\mathbf{N}_j = N(y_j^{[n]})$ \\
        \hspace{1.5em} \text{{\bf for} j = 1 to $q$}  {\color{gray} \# {\bf parallelizable loop}} \\  
       	\hspace{3em} $y_j^{[n+1]} = \varphi_0(r\mathbf{L}\eta_j(0))y_1^{[n+1]} + \sum_{k=0}^{q-1} (\eta_j(0))^{k+1} \varphi_{k+1}(r\mathbf{L}\eta_j(0))  L_N^{(k)}(-1) $ \\[0.5em] \hline

\end{tabularx}
\vspace{0.1em}
\caption{Pseudocode for a single timestep of a partitioned composite EPBM method (\ref{eq:leg_method_composite}) where the method $M(\alpha)$ is constructed using the parameters from Table \ref{tab:leg_methods}. All for loops that can be run in parallel have been explicitly annotated.}
	\label{tab:EPBM_pseudocode}
\end{table}

		\subsubsection{Parallel Adams EPBMs with imaginary equispaced nodes}
		\label{subsubsection:imaginary-nodes}
		
		We demonstrate the flexibility of the polynomial framework by reusing the PMFC strategy to construct two example parallel EPBM methods using the set of $m$ imaginary equispaced nodes
			\begin{align}
				\chi_{j,m} = i\left(1 - \frac{2(j-1)}{m-1}\right), && j = 1, \ldots, k.
				\label{eq:imag_equi_nodes}
			\end{align}
			The parameters for the two methods are described in Table \ref{tab:iequi_methods}. The method A is identical to the method from Table \ref{tab:leg_methods} except that the Legendre nodes have been replaced with imaginary equispaced nodes, while method $B$ is the exponential equivalent of the BAM method from \cite{buvoli2019constructing}. Serial imaginary methods can also be constructed by using the SMFC and SMVC strategies for imaginary nodes that are proposed in \cite{buvoli2018polynomial}.

		\begin{table}[h!]
			\renewcommand{\arraystretch}{1.5}
			\begin{tabular}{rccc}
															& {\bf method A} 										& \hspace{1em} & {\bf method B} \\
				$\{z_j\}_{j=1}^q$ : 							& 	$\{-1, \chi_{1,q-1}, \ldots, \chi_{q-1,q-1}\}$ 	& 	& $\{ \chi_{j,q} \}_{j=1}^q$ \\
				$L^{[j]}_y(\tau)$ :							& 	$L_y(\tau)$ from (\ref{eq:lagrange_Ly_formula}) & 	&  $L_y(\tau)$ from (\ref{eq:lagrange_Ly_formula}) \\
				$L^{[j]}_N(\tau)$ or $L^{[j]}_R(\tau)$ :		&	PMFC$_2$										& 	& PMFC$_1$ \\
				$\{b_j\}_{j=1}^q$ : 							& 	$\{z_1\}_{j=1}^q$ 								& 	& $\{z_j\}_{j=1}^q$					
			\end{tabular}
			\caption{Parameters for EPBMs with imaginary equispaced nodes. The constants $\chi_{i,j}$ represent the imaginary equispaced formula from (\ref{eq:imag_equi_nodes}).} 
			\label{tab:iequi_methods}
		\end{table}

		As we have shown here, constructing methods with complex-valued nodes is no more difficult than constructing methods with real-valued nodes. However, extending into the complex plane introduces complexities including: (1) the solution must be sufficiently analytic off the real time-line, (2) matrix exponentials must also be efficient to compute at complex times, and (3) the code for the ODE right-hand-side must accept complex arguments.
		
		For diagonally implicit PBMs from \cite{buvoli2019constructing}, extending into the complex plane provided significantly improved stability compared to real-valued nodes. However, for exponential PBMs we can obtain good stability with real valued nodes while avoiding the aforementioned issues. With the exception of a single numerical experiment contained in the appendix, we therefore focus entirely on EPBMs with real-valued nodes.

\section{Linear stability}
\label{sec:stability_accuracy} 

We now briefly analyze the linear stability properties of partitioned and unpartitioned Legendre EPBMs from Table \ref{tab:leg_methods}. Like all unpartitioned integrators, any unpartitioned EPBM will integrate the classical Dalhquist test problem $y' = \lambda y$ exactly, and therefore its linear stability region is always equal to the left-half-plane $\text{Re}(h\lambda) \le 0$. To analyze linear stability for partitioned methods we use the generalized Dahlquist test problem \cite{ascher1995implicit, izzo2017highly, jackiewicz2017construction, cardone2015construction, cox2002ETDRK4, krogstad2005IF, grooms2011IMEXETDCOMP, sandu2015generalized}
    \begin{align}
    	y' = \lambda_1 y + \lambda_2 y,
    	\label{eq:partitioned_dahlquist_test_problem}
    \end{align}
	where the terms $\lambda_1 y$ and $\lambda_2 y$ are meant to respectively represent the linear and nonlinear terms. When  applied to  (\ref{eq:partitioned_dahlquist_test_problem}), an EPBM with $q$ inputs reduces to the matrix iteration
	\begin{align}
		\mathbf{y}_{n+1} = \mathbf{M}(z_1, z_2, \alpha) \mathbf{y}_n	
	\end{align}
	where $\mathbf{M}$ is a $q\times q$ matrix, $z_1 = h \lambda_1$, $z_2 = h \lambda_2$, and $h$ denotes the timestep. The stability region $S$ is the subset of $\mathbb{C}^2 \times \mathbb{R}^+$ where $\mathbf{M}(z_1, z_2, \alpha)$ is power bounded, so that 
		\begin{align}
			S = \left\{z ~:~ \sup_{n\in \mathbb{N}} \| \mathbf{M}(z_1, z_2, \alpha)^n \| < \infty	 \right\}.
		\end{align}
		The matrix $\mathbf{M}(z_1, z_2, \alpha)$ is power bounded if its eigenvalues lie inside the closed unit disk, and any eigenvalues of magnitude one are non-defective. The stability region $\mathcal{S}$ is five-dimensional and cannot be shown directly. Instead, we overlay two-dimensional slices of the stability regions formed by fixing $z_1$ and $\alpha$ while allowing $z_2$ to vary. The choices for $z_1$ are:
	\begin{enumerate}
		\item $z_1 \in -1 \cdot [0,30]$: negative real-values mimic a problem with linear dissipation.
		\item $z_1 \in i \cdot [0,30]$: imaginary values mimic linear advection or dispersion.
		\item $z_1 \in \exp\left(\frac{3\pi i}{4}\right) \cdot [0,30]$: complex-values mimic both dispersion and dissipation.
	\end{enumerate}

We compare the linear stability regions of Legendre EPBMs to those of exponential Adams-Bashforth (EAB)  \cite{beylkin1998ELP} and the fourth-order exponential Runge-Kutta (ERK) method from \cite{cox2002ETDRK4}. It is interesting to compare EPBMs to EABs since both methods are constructed using high-order Lagrange polynomials. On the other hand, ERK methods serve as a good benchmark for gauging stability since they are stable on a variety of problems \cite{montanelli2016solving, KassamTrefethen05ETDRK4, cox2002ETDRK4, grooms2011IMEXETDCOMP}. 

\begin{figure}[hbt!]

	\begin{tabular}{lccc}
		& {\bf Dissipative} & {\bf Mixed}  & {\bf Oscillatory} \\[0.5em]
		& {$z_1 = -r$} & {$z_1 = r\exp(\frac{3 \pi i}{4})$}  & {$z_1 = ir$} \\[1em]
		\rotatebox{90}{\hspace{2em} EAB4} & 
		\includegraphics[width=0.28\linewidth, valign=b]{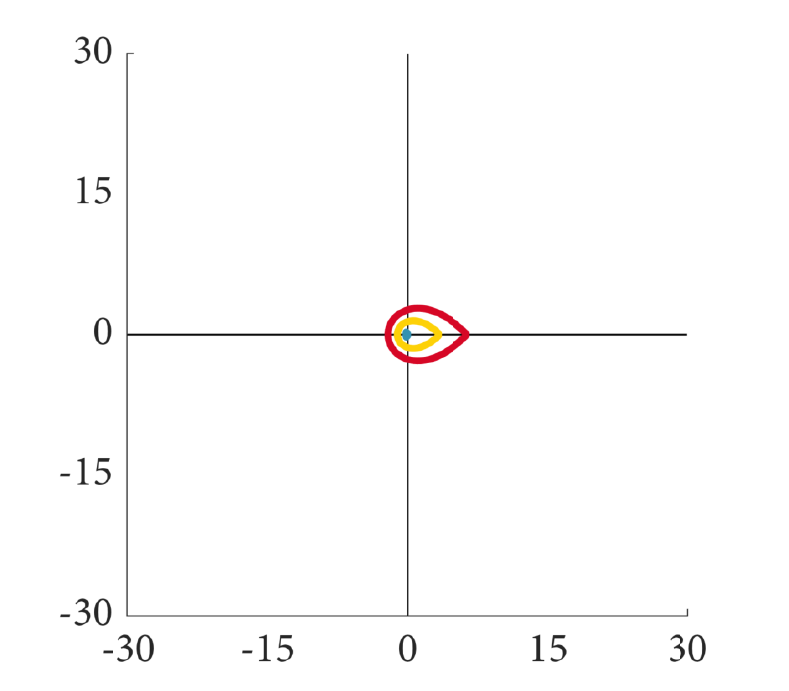} &
		\includegraphics[width=0.28\linewidth, valign=b]{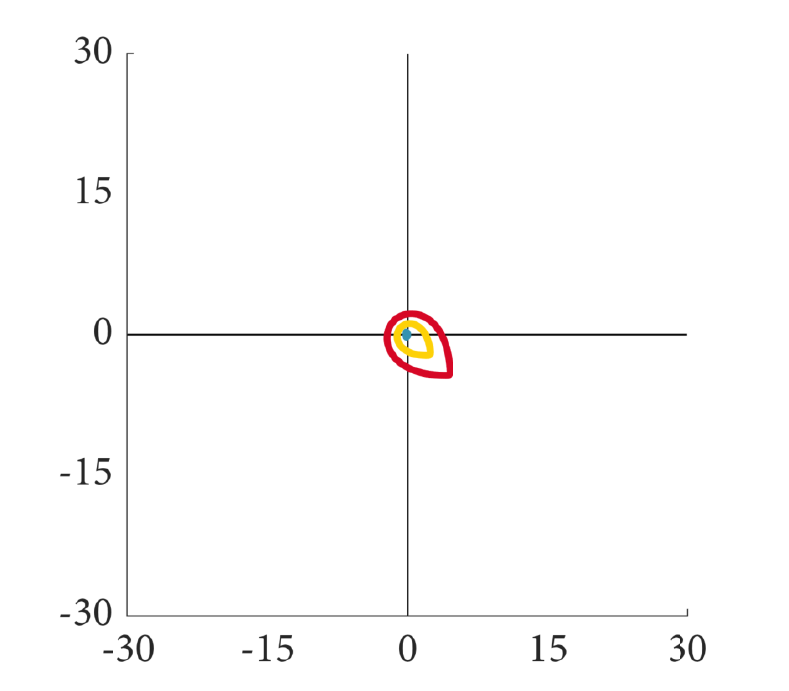} &
		\includegraphics[width=0.28\linewidth, valign=b]{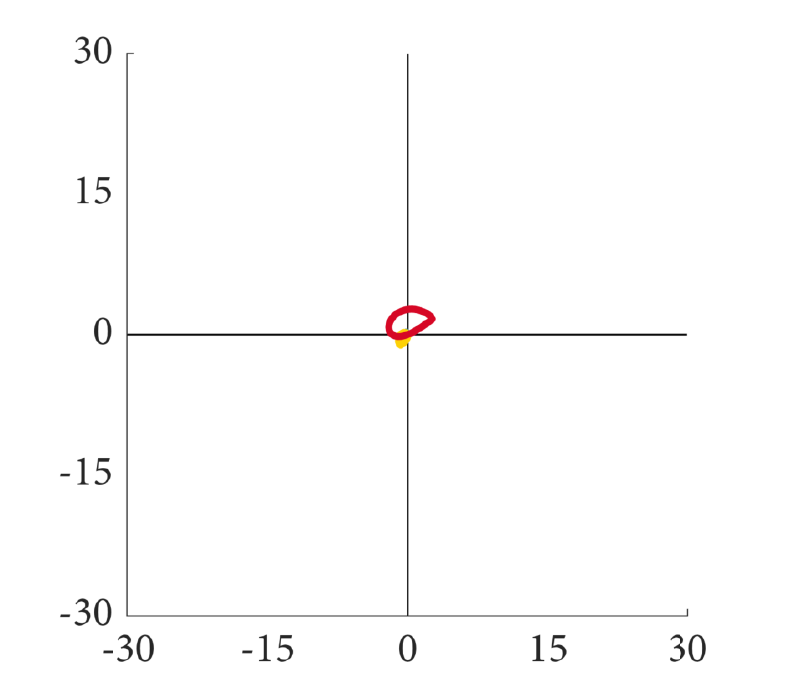} \\
		\rotatebox{90}{\hspace{2em} ERK4} & 
		\includegraphics[width=0.28\linewidth, valign=b]{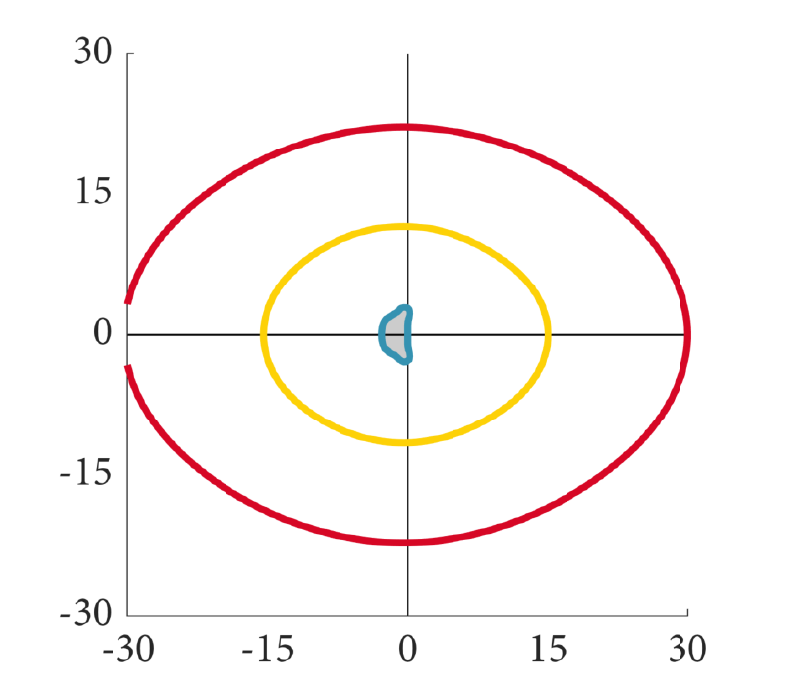} &
		\includegraphics[width=0.28\linewidth, valign=b]{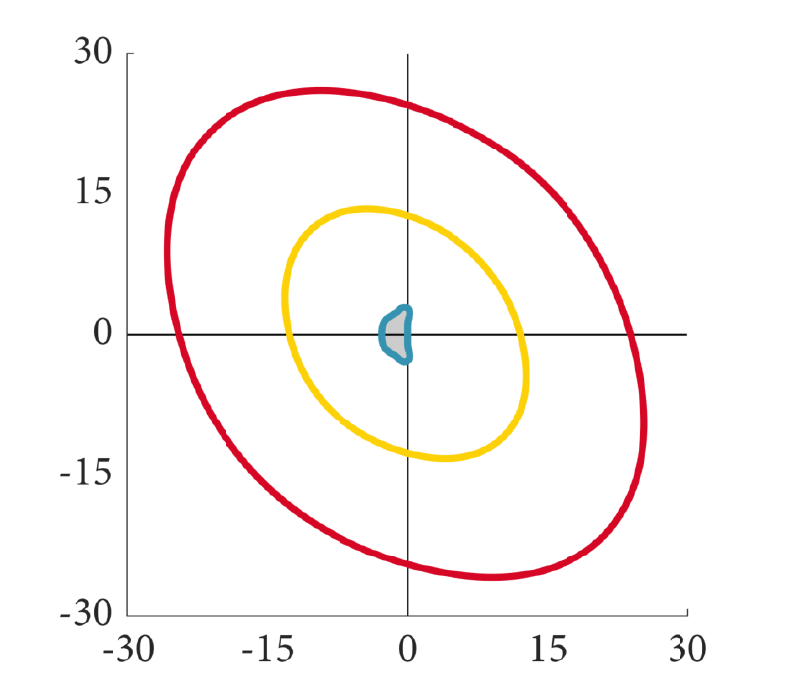} &
		\includegraphics[width=0.28\linewidth, valign=b]{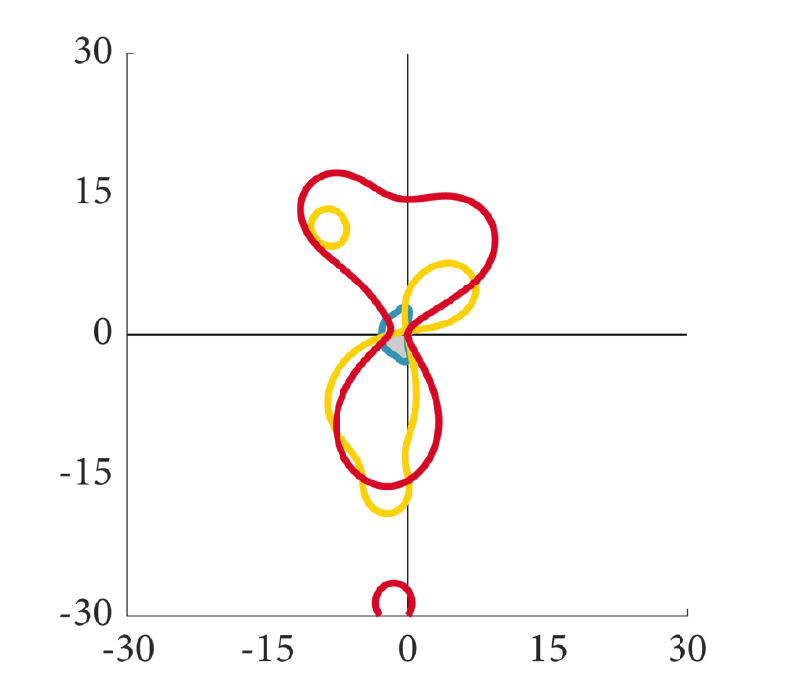} \\
		\rotatebox{90}{\hspace{2em} EPBM4} & 
		\includegraphics[width=0.28\linewidth, valign=b]{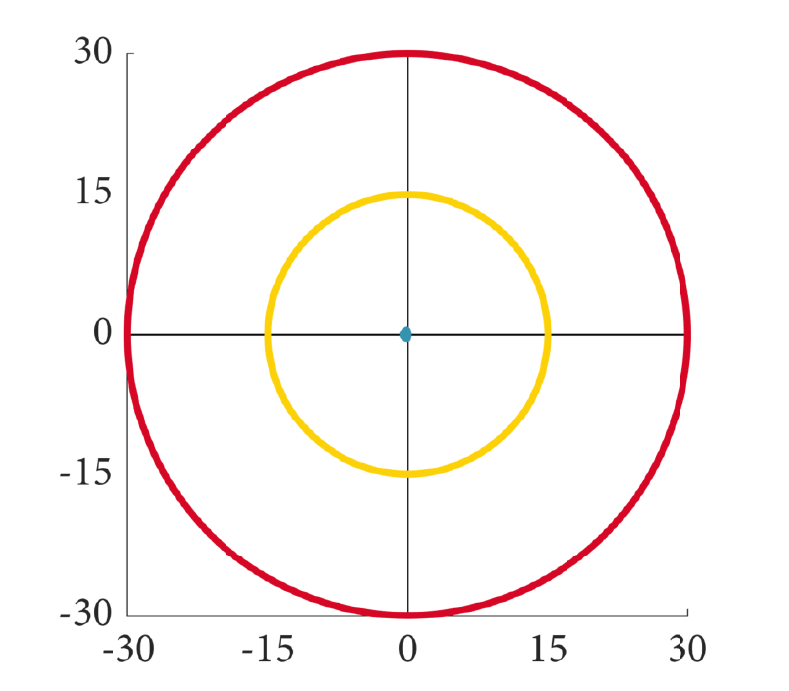} &
		\includegraphics[width=0.28\linewidth, valign=b]{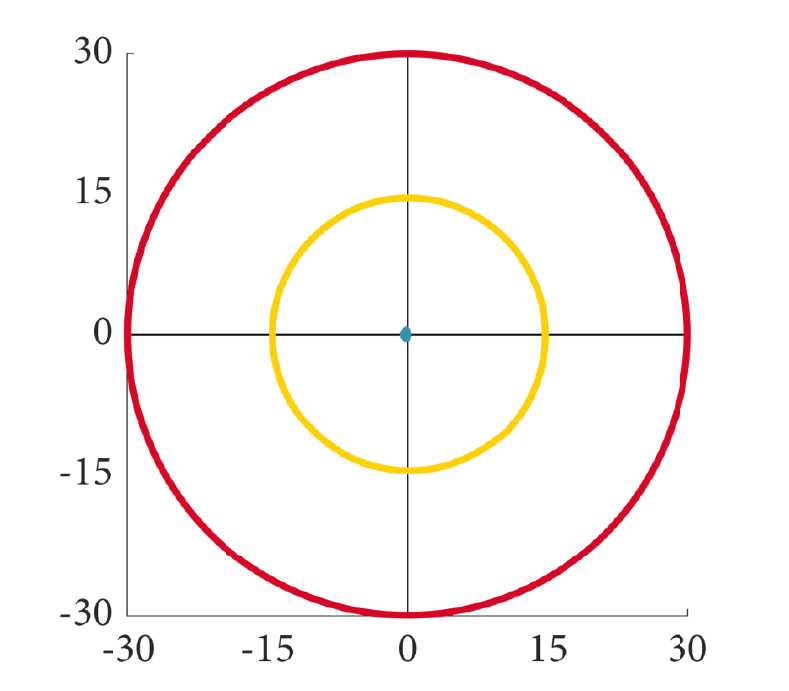} &
		\includegraphics[width=0.28\linewidth, valign=b]{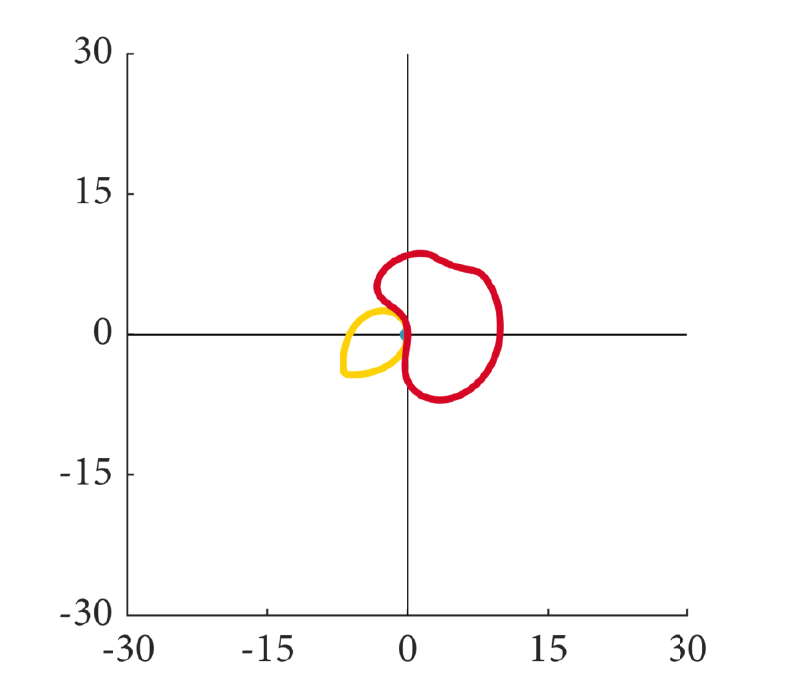} \\
		\rotatebox{90}{\hspace{2.5em} EPBM8} & 
		\includegraphics[width=0.28\linewidth, valign=b]{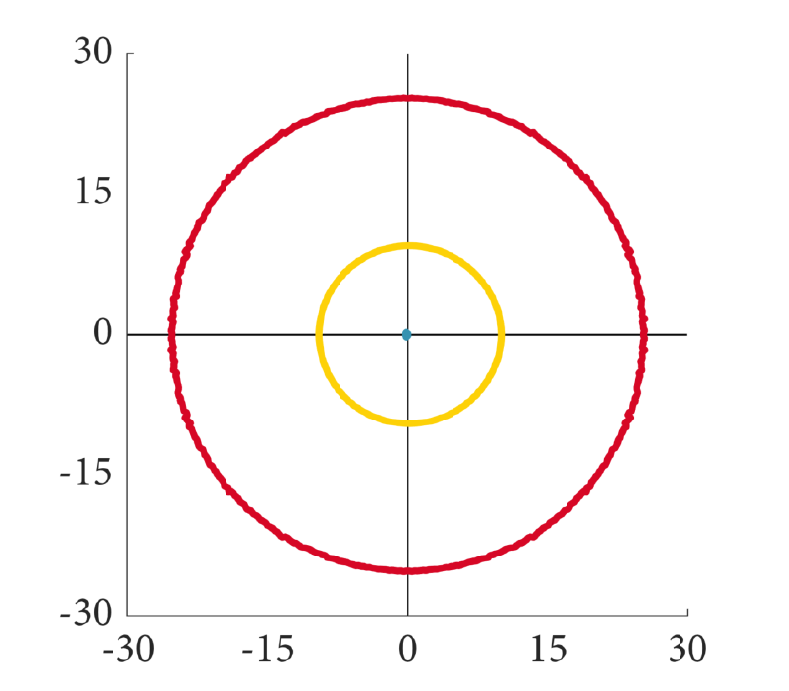} &
		\includegraphics[width=0.28\linewidth, valign=b]{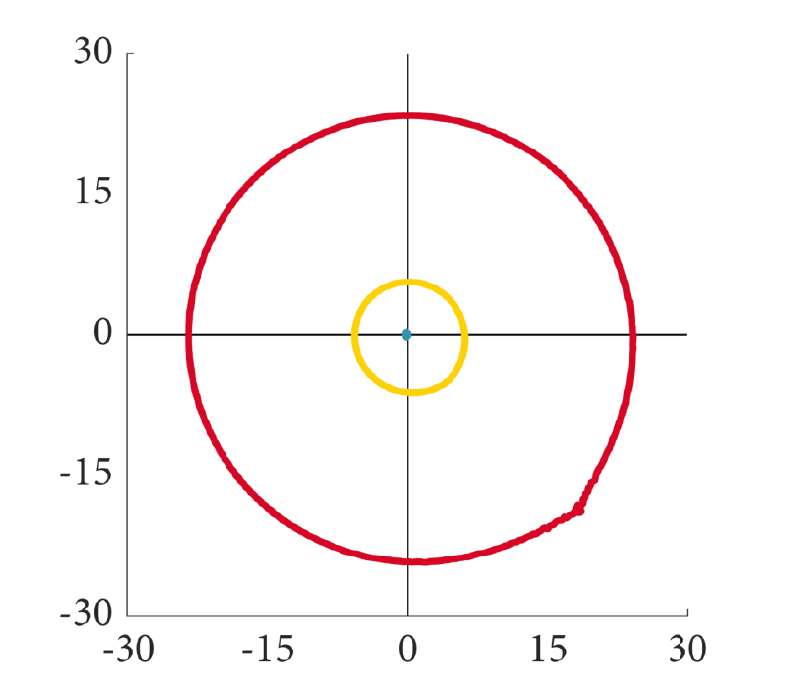} &
		\includegraphics[width=0.28\linewidth, valign=b]{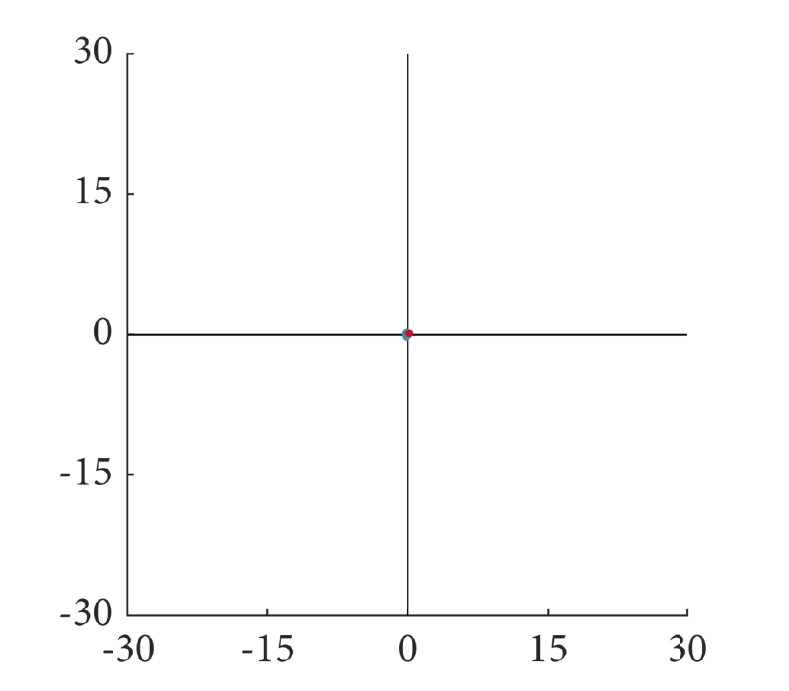}
	\end{tabular}
	
	\vspace{1em}
	\begin{center}
	\begin{tabular}{c}
	{\tiny \textcolor{plot_blue}{\hdashrule[0.2ex]{2em}{2pt}{}} $r = 0$ \hspace{1em}}
	{\tiny \textcolor{plot_yellow}{\hdashrule[0.2ex]{2em}{2pt}{}} $r = 15$ \hspace{1em}}
	{\tiny \textcolor{plot_red}{\hdashrule[0.2ex]{2em}{2pt}{}} $r = 30$} \hspace{1em}
	\end{tabular}
	\end{center}
	\vspace{1em}

	\caption{Linear stability regions of the fourth-order EAB \cite{beylkin1998ELP}, the fourth-order ERK4 \cite{cox2002ETDRK4} and the fourth-order and eighth-order Legendre EPBMs from Table \ref{tab:leg_methods}. Each row corresponds to a method and each column to one of the three choices for $z_1$ described in Section \ref{sec:stability_accuracy}. The colored contours represent the stability regions for $z_1$ with magnitudes of $0$, $15$, and $30$. We did not include the stability regions for eighth-order exponential Adams-Bashforth in this plot since they are too small to be visible at this scale.}
	\label{fig:stability_large}
	
\end{figure}

\begin{figure}[hbt!]
	\begin{tabular}{lccc}
		& {\bf Dissipative} & {\bf Mixed}  & {\bf Oscillatory} \\
		& {$z_1 = -r$} & {$z_1 = r\exp(\frac{3 \pi i}{4})$}  & {$z_1 = ir$} \\
		\rotatebox{90}{\hspace{2em} EAB4} & 
		\includegraphics[width=0.28\linewidth, valign=b]{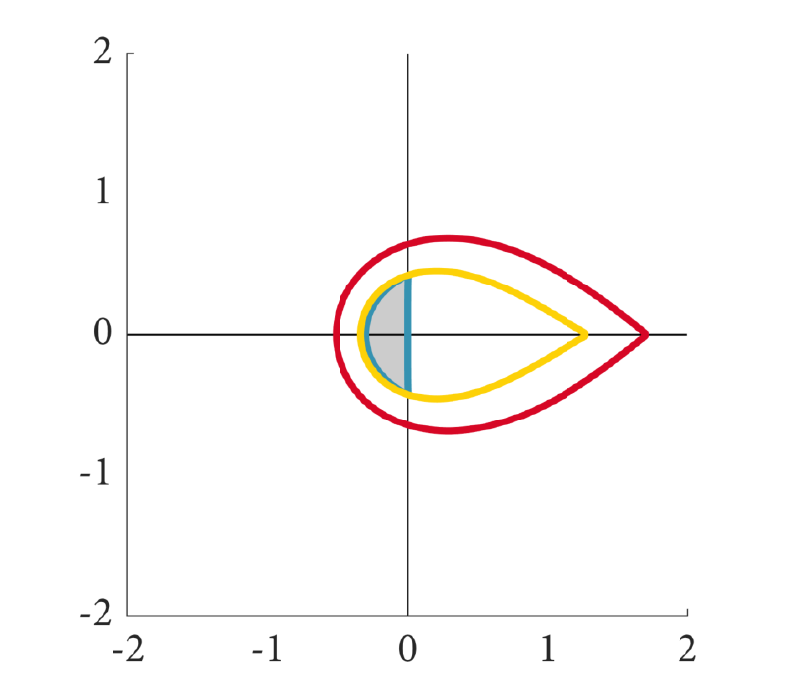} &
		\includegraphics[width=0.28\linewidth, valign=b]{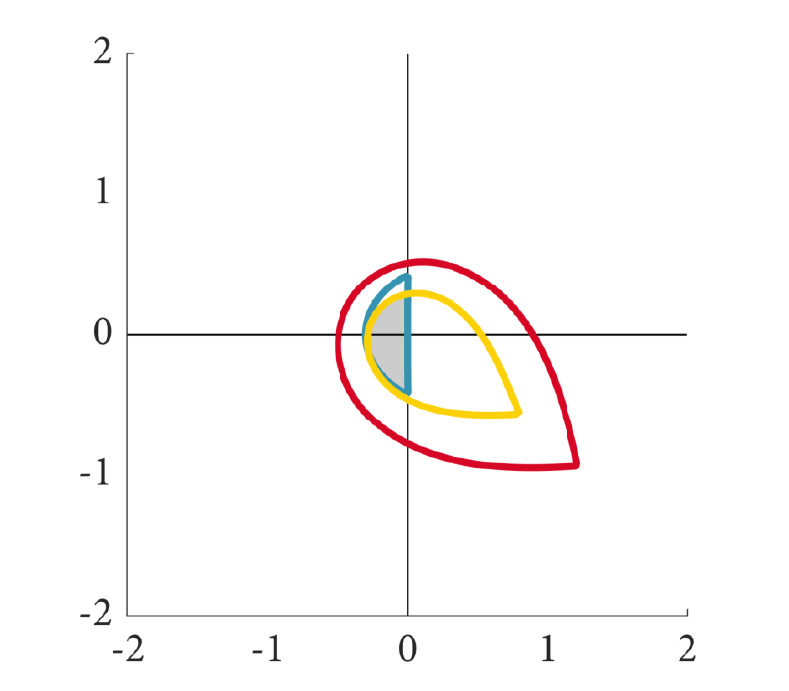} &
		\includegraphics[width=0.28\linewidth, valign=b]{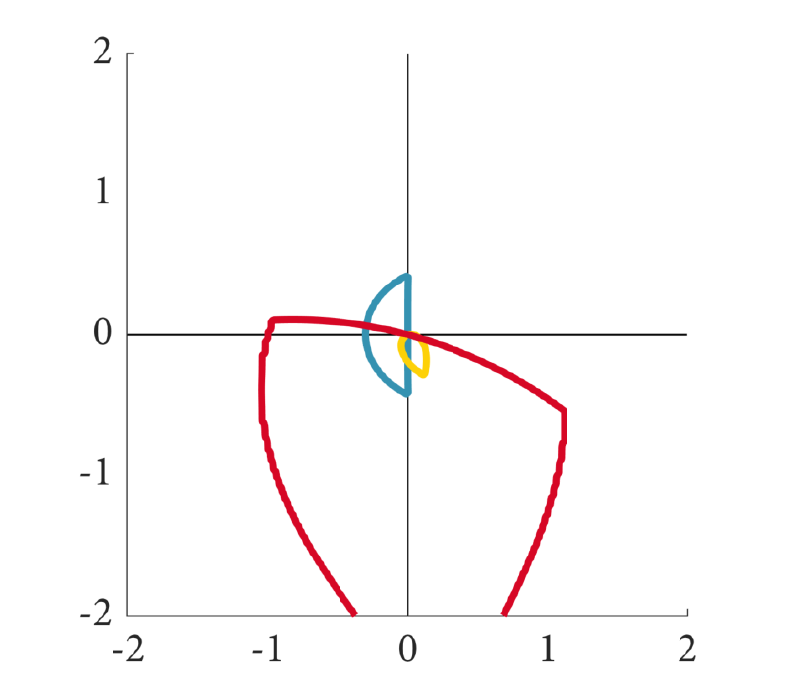} \\
		\rotatebox{90}{\hspace{2.5em} EAB8} & 
		\includegraphics[width=0.28\linewidth, valign=b]{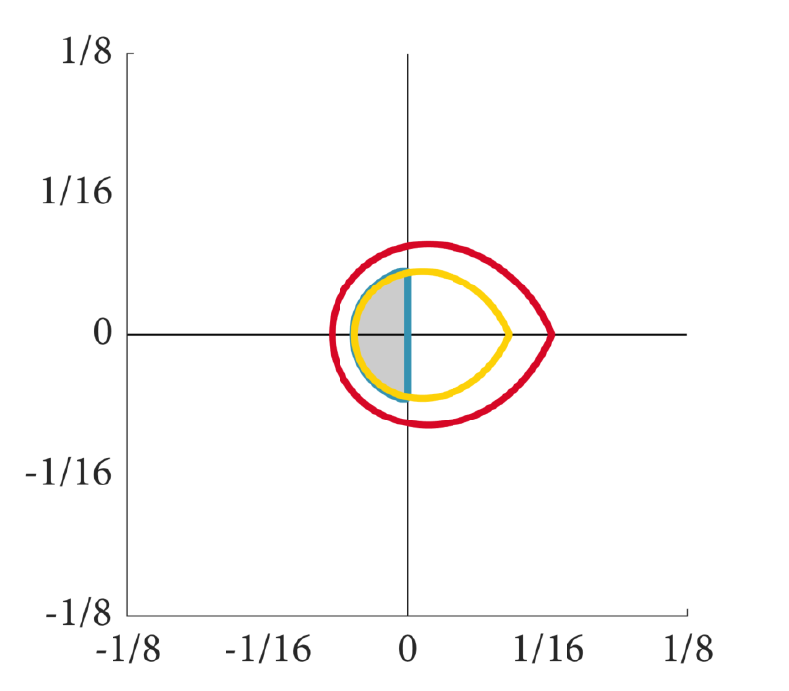} &
		\includegraphics[width=0.28\linewidth, valign=b]{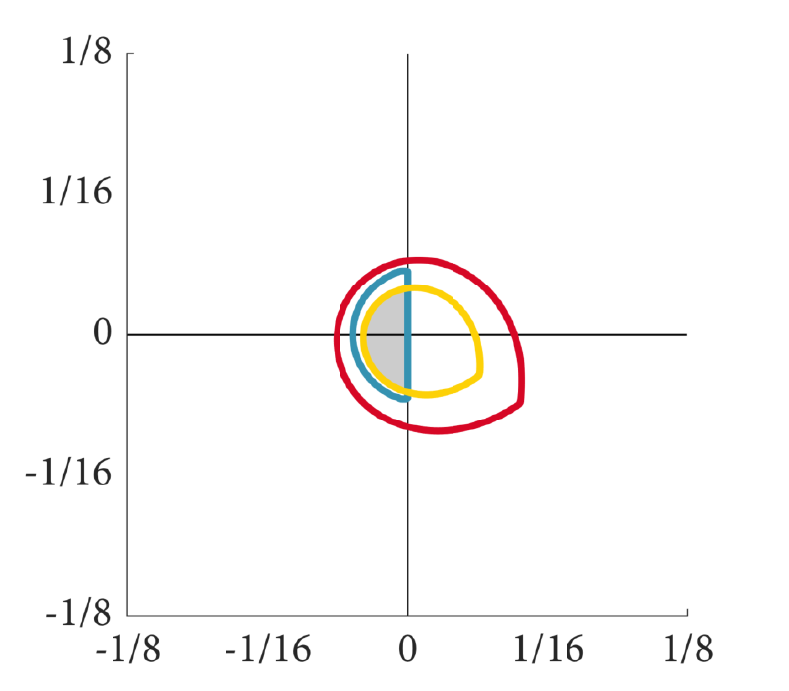} &
		\includegraphics[width=0.28\linewidth, valign=b]{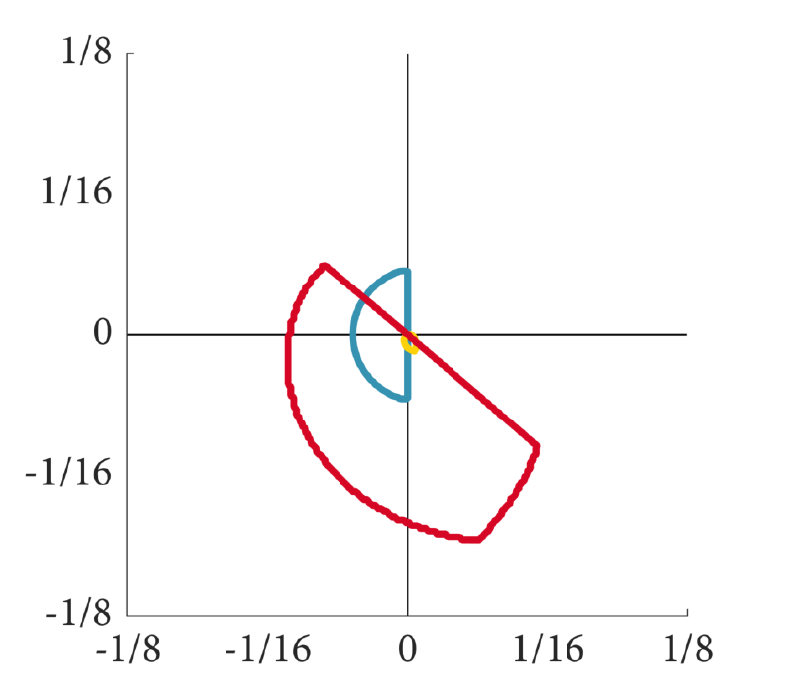} \\
		\rotatebox{90}{\hspace{2em} ERK4} & 
		\includegraphics[width=0.28\linewidth, valign=b]{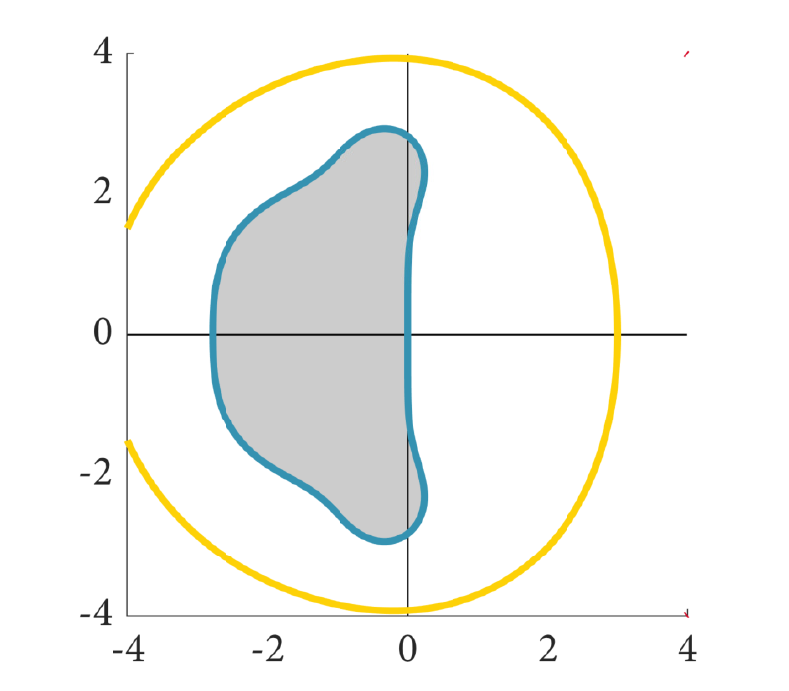} &
		\includegraphics[width=0.28\linewidth, valign=b]{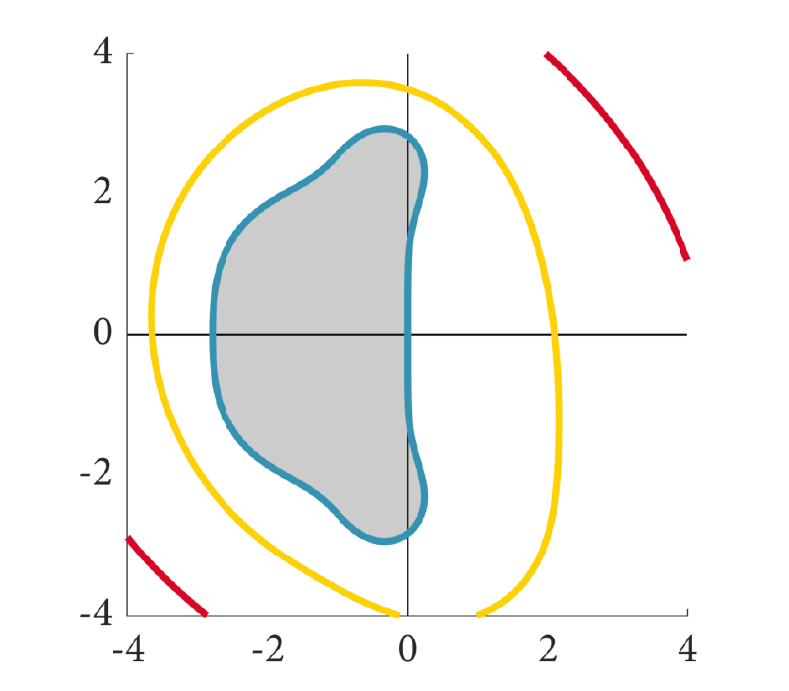} &
		\includegraphics[width=0.28\linewidth, valign=b]{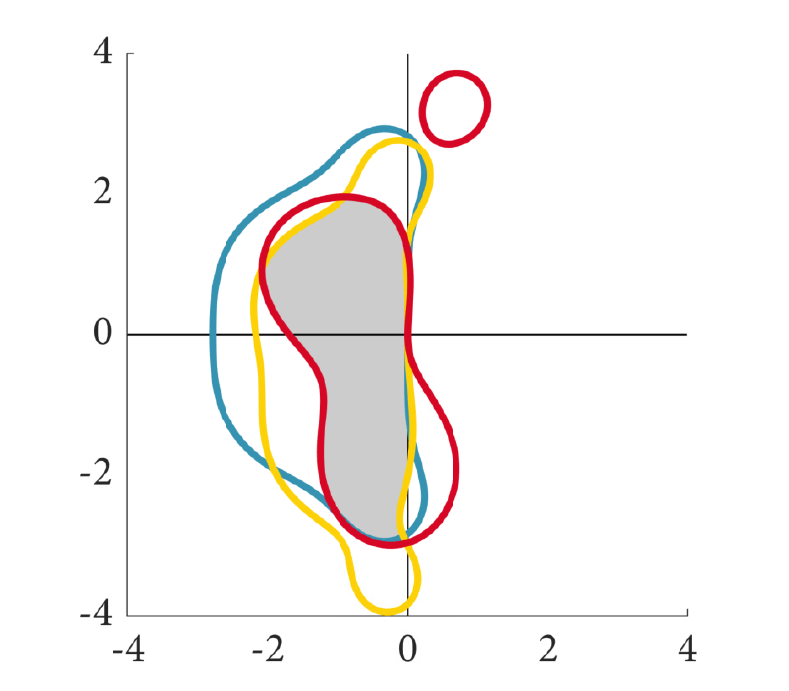} \\
		\rotatebox{90}{\hspace{2em} EPBM4} & 
		\includegraphics[width=0.28\linewidth, valign=b]{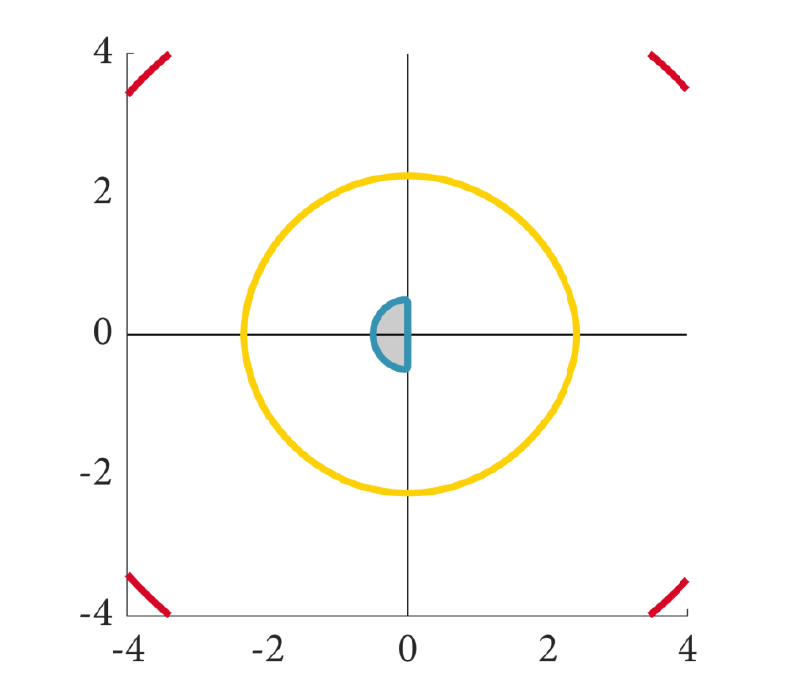} &
		\includegraphics[width=0.28\linewidth, valign=b]{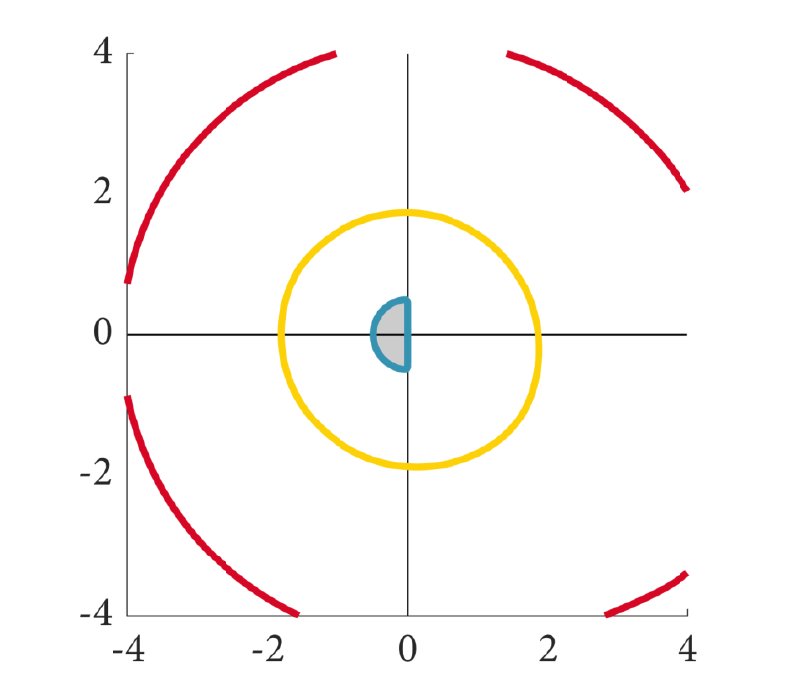} &
		\includegraphics[width=0.28\linewidth, valign=b]{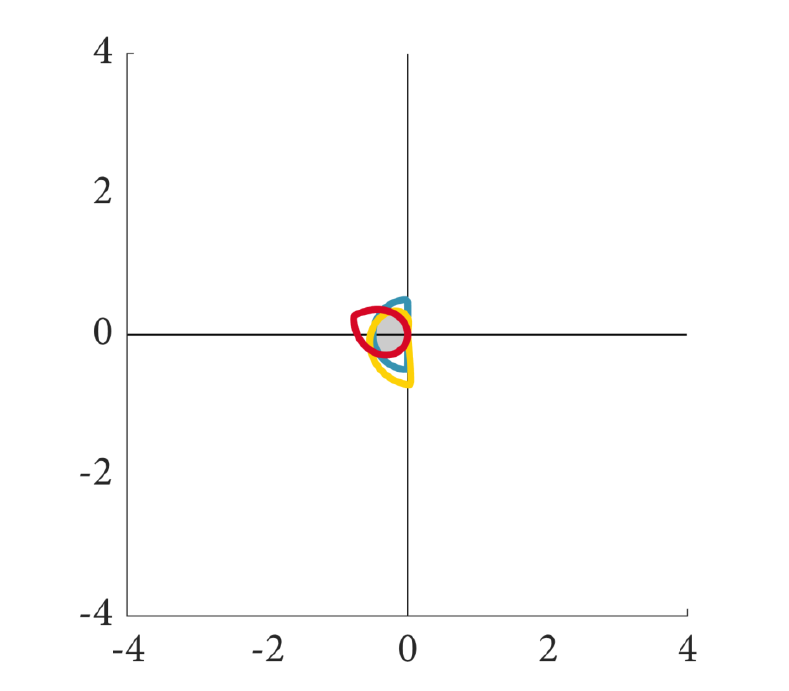} \\
		\rotatebox{90}{\hspace{2.5em} EPBM8} & 
		\includegraphics[width=0.28\linewidth, valign=b]{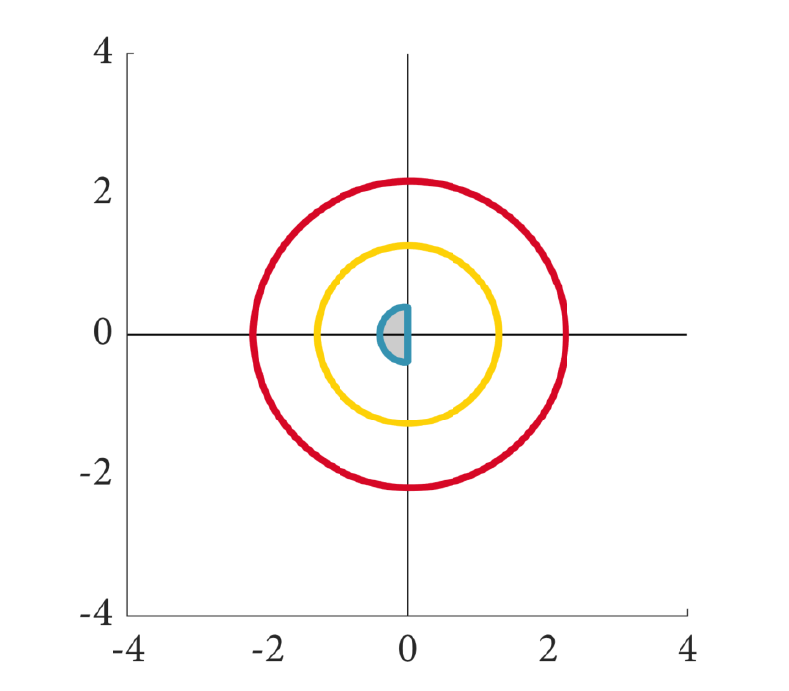} &
		\includegraphics[width=0.28\linewidth, valign=b]{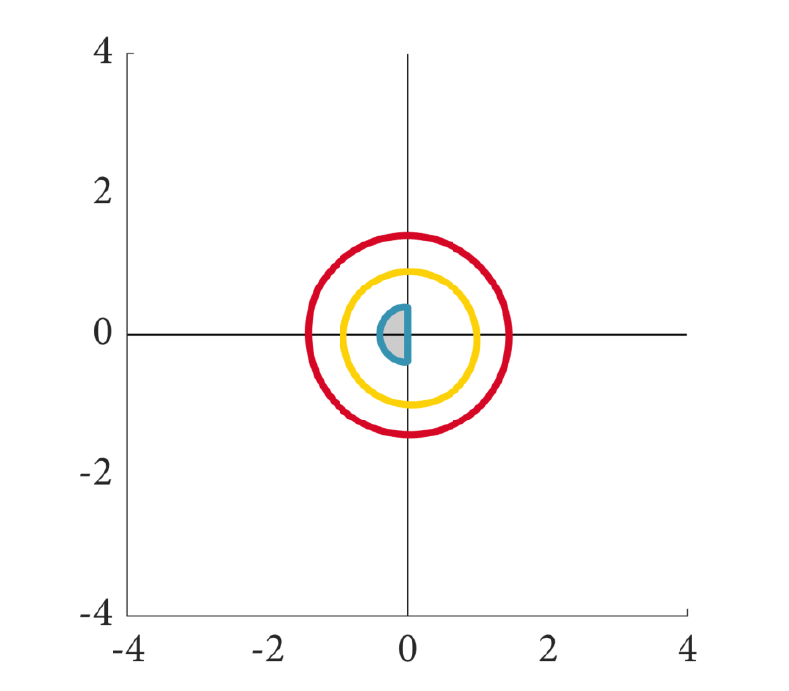} &
		\includegraphics[width=0.28\linewidth, valign=b]{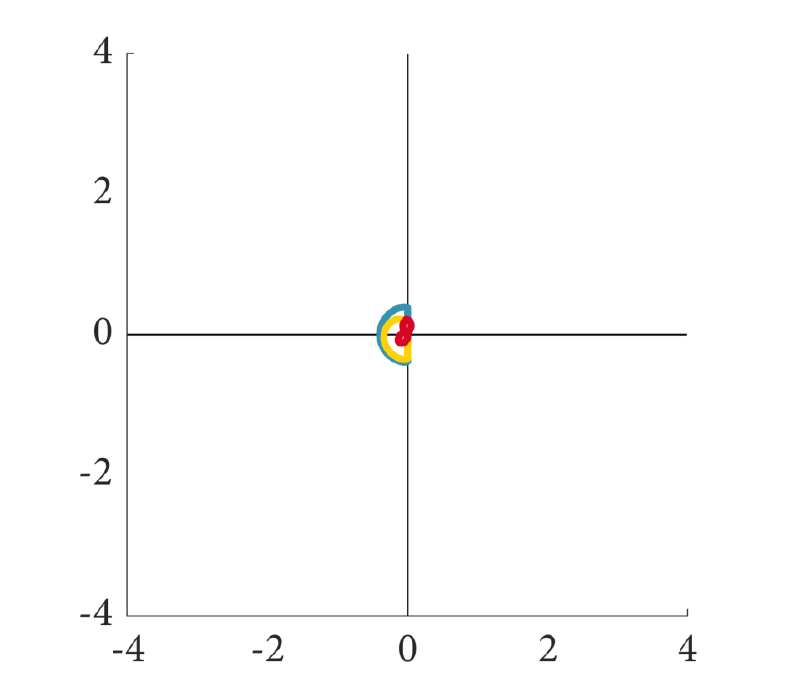}
	\end{tabular}
	
	\vspace{1em}
	\begin{center}
		\begin{tabular}{c}
			{\tiny \textcolor{plot_blue}{\hdashrule[0.2ex]{2em}{2pt}{}} $\left|z_1\right| = 0$ \hspace{1em}}
			{\tiny \textcolor{plot_yellow}{\hdashrule[0.2ex]{2em}{2pt}{}} $\left|z_1\right| = 3$ \hspace{1em}}
			{\tiny \textcolor{plot_red}{\hdashrule[0.2ex]{2em}{2pt}{}} $\left|z_1\right| = 6$} \hspace{1em}
		\end{tabular}
	\end{center}

	\caption{Magnified stability regions from Figure \ref{fig:stability_large}. The axis used in the plots for the EAB4 and EAB8 are magnified by a factor of 2 and 32, respectively, compared to those of the ERK4 and the Legendre EPBMs. The colored contours now denote stability regions where $z_1$ has magnitude $0$, $3$, and $6$.}
	\label{fig:stability_small}
\end{figure}

In Figure \ref{fig:stability_large} we show the stability regions for fourth and eighth order Legendre EPBMs with $\alpha = 2$ compared against those of ERK4 and the fourth and eighth order EAB methods. In Figure \ref{fig:stability_small} we also show magnified plots that highlight the stability properties near the origin.

When dissipation is present, the stability regions for EPBMs are significantly larger than those of EAB schemes. The difference between the two methods increases with order, with the eighth-order EPBM possessing a stability region that is approximately thirty-two times larger than the eighth-order EAB method for small $|z_1|$. For large values of $|z_1|$, the stability regions of EPBMs are comparable to those of the fourth-order Runge-Kutta, however, for smaller $|z_1|$, Runge-Kutta integrators display larger stability regions.

\begin{figure}[hbt!]
	\begin{tabular}{lccc}
		& { $\kappa = 0$} & { $\kappa = 1$}  & {$\kappa = 2$} \\
		\rotatebox{90}{\hspace{2em} EPBM4} & 
		\includegraphics[width=0.28\linewidth, valign=b]{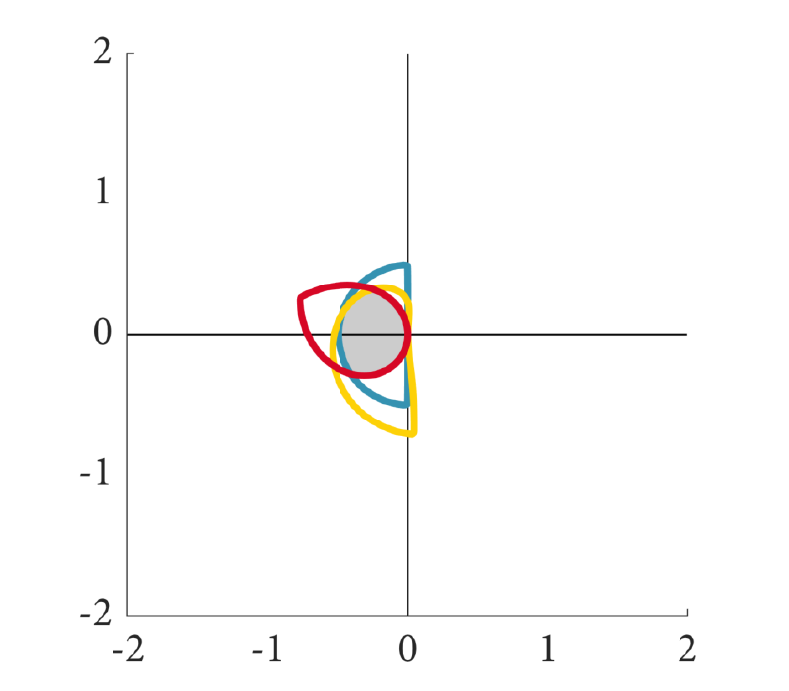} &
		\includegraphics[width=0.28\linewidth, valign=b]{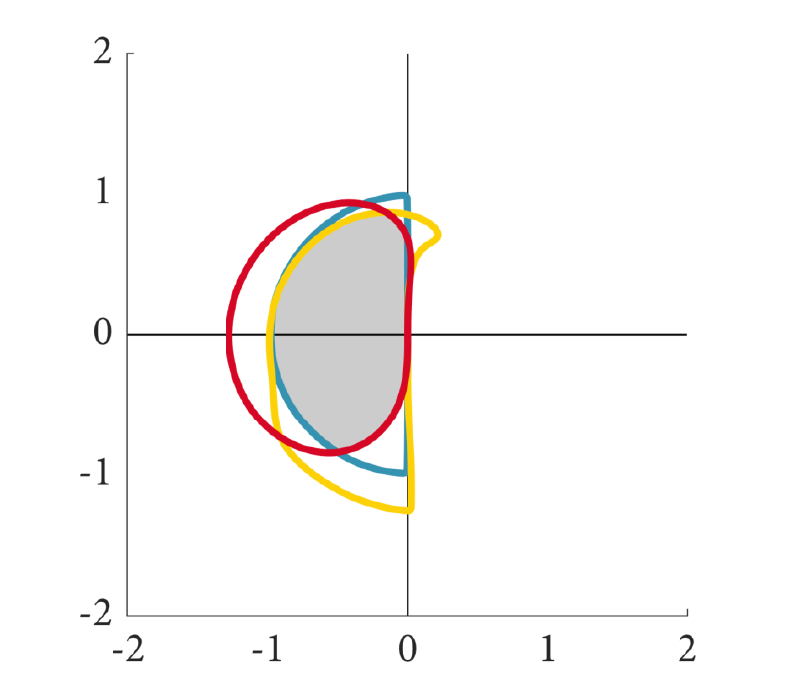} &
		\includegraphics[width=0.28\linewidth, valign=b]{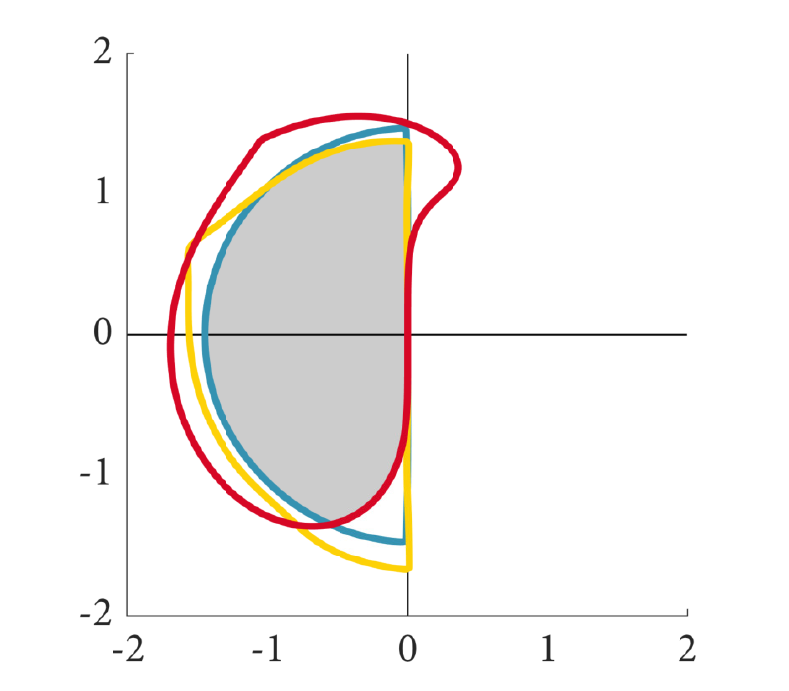} \\
		\rotatebox{90}{\hspace{2.5em} EPBM6} & 
		\includegraphics[width=0.28\linewidth, valign=b]{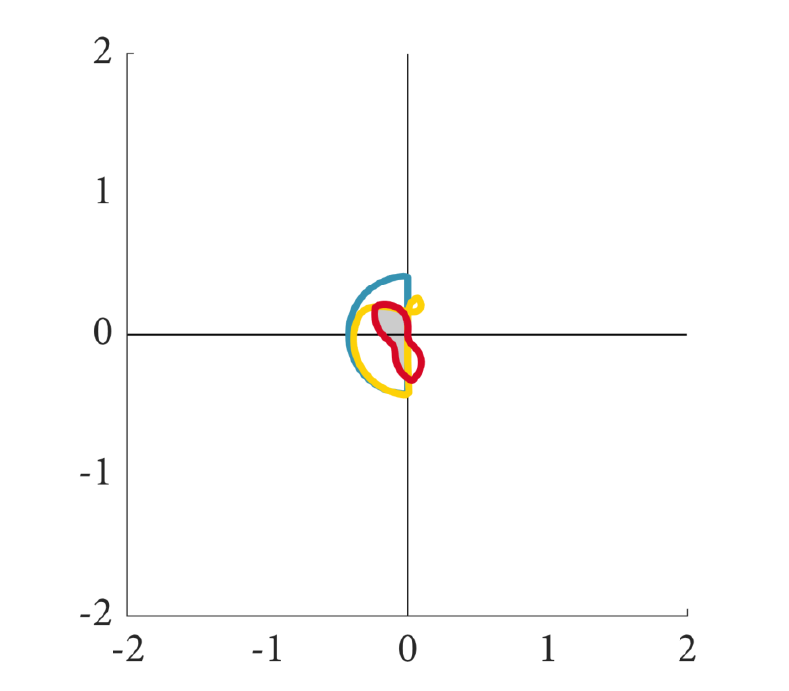} &
		\includegraphics[width=0.28\linewidth, valign=b]{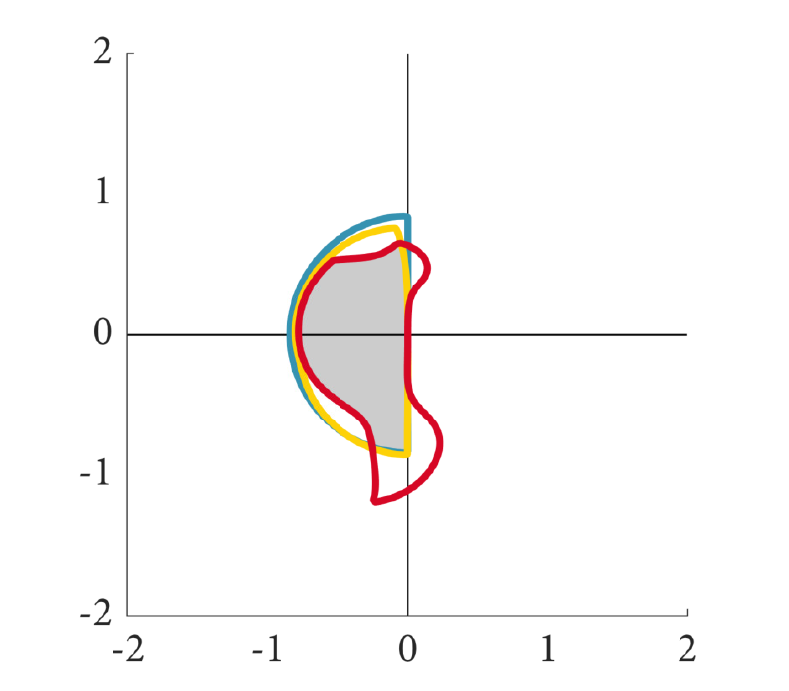} &
		\includegraphics[width=0.28\linewidth, valign=b]{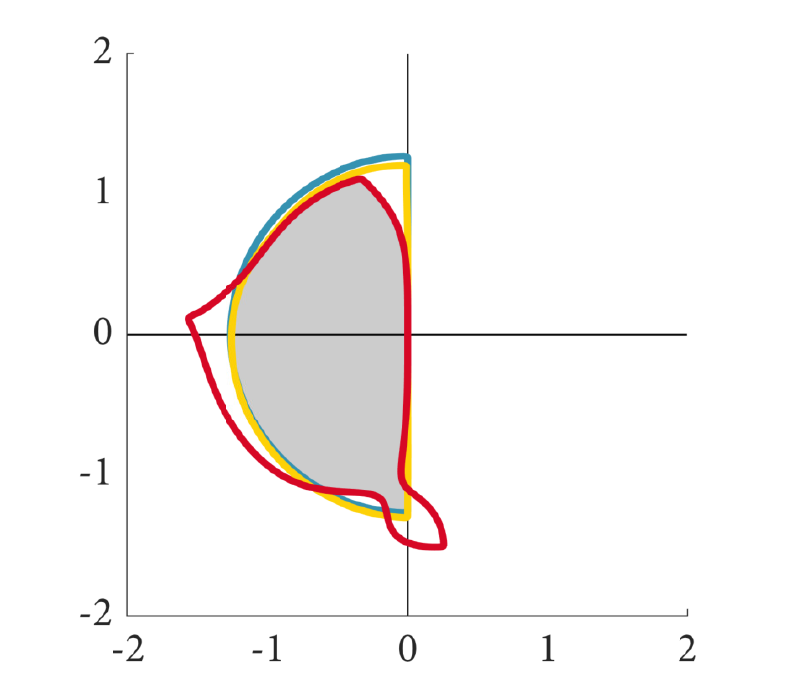} \\
		\rotatebox{90}{\hspace{2.5em} EPBM8} & 
		\includegraphics[width=0.28\linewidth, valign=b]{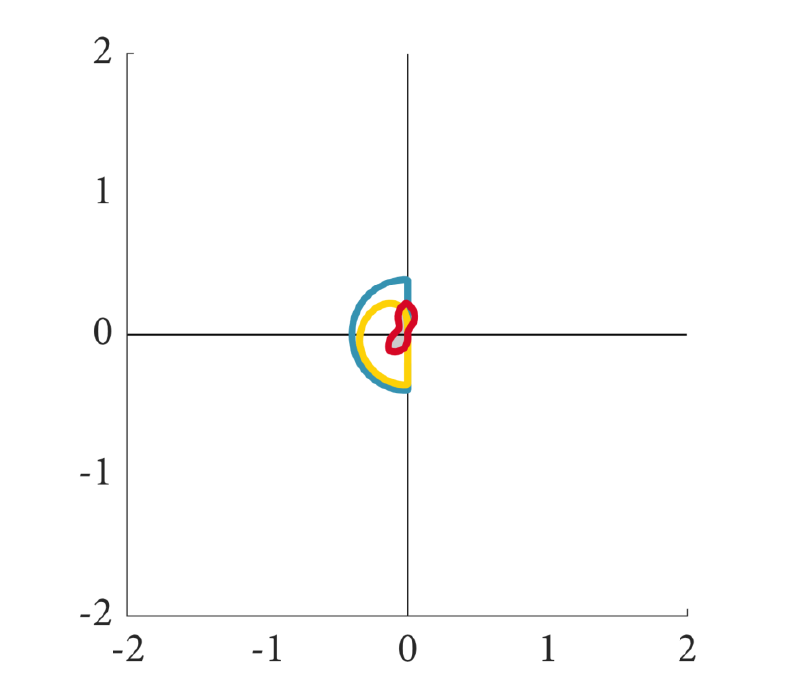} &
		\includegraphics[width=0.28\linewidth, valign=b]{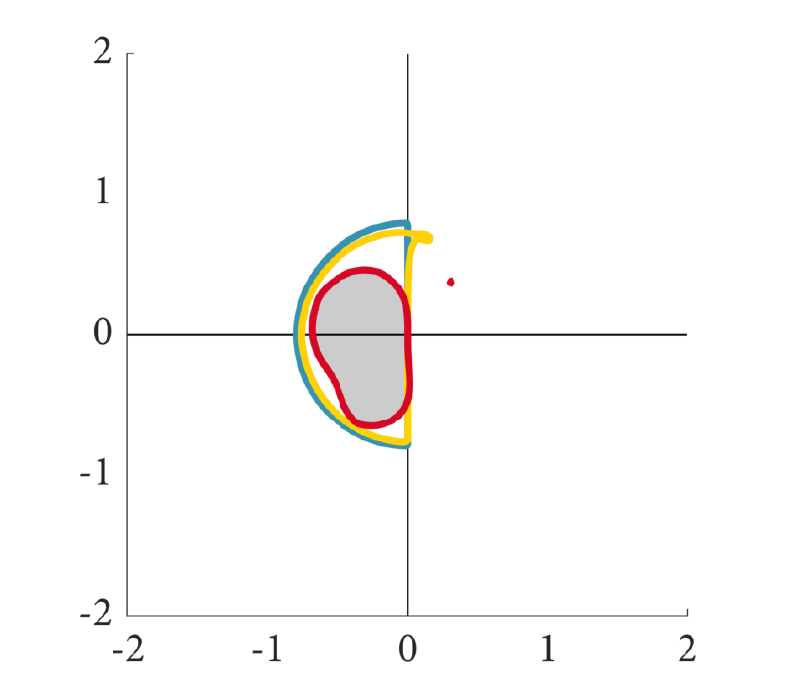} &
		\includegraphics[width=0.28\linewidth, valign=b]{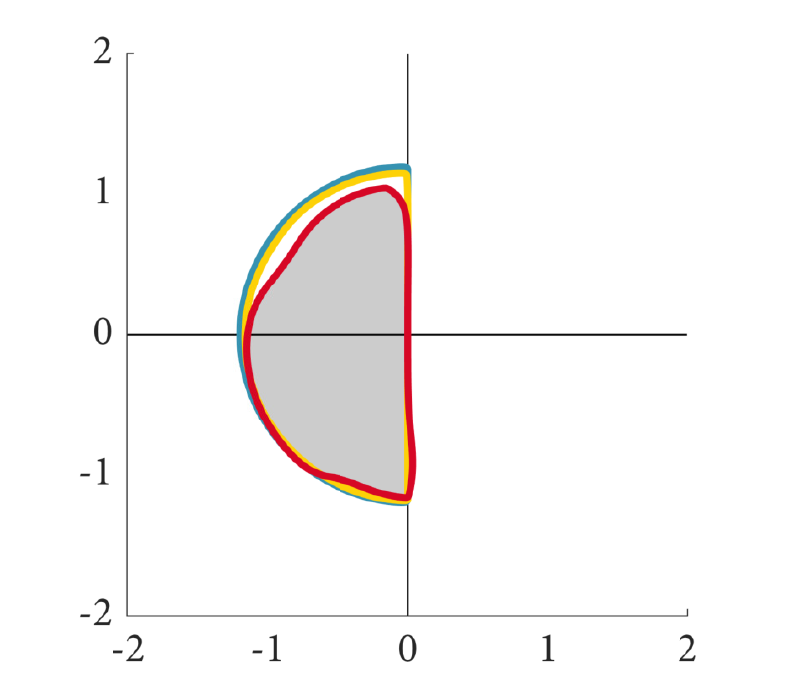} \\
	\end{tabular}
	
	\vspace{1em}
		\begin{center}
			\begin{tabular}{c}
			{\tiny \textcolor{plot_blue}{\hdashrule[0.2ex]{2em}{2pt}{}} $\left|z_1\right| = 0$ \hspace{1em}}
			{\tiny \textcolor{plot_yellow}{\hdashrule[0.2ex]{2em}{2pt}{}} $\left|z_1\right| = 3$ \hspace{1em}}
			{\tiny \textcolor{plot_red}{\hdashrule[0.2ex]{2em}{2pt}{}} $\left|z_1\right| = 6$} \hspace{1em}
			\end{tabular}
		\end{center}
	\vspace{1em}

	\caption{Magnified stability regions for composite Legendre methods (\ref{eq:leg_method_composite}) on the oscillatory Dahlquist test problem where $z_1 = ir$. The composite methods takes one propagation step with $\alpha = 2$ and $k$ iteration steps with $\alpha = 0$. The stability regions increase substantially even after a single iteration step.}
	\label{fig:stability_iteration}	
\end{figure}

For the purely oscillatory Dahlquist test problem, all three families of integrators have small stability regions that contract as $|z_1|$ increases. For small $|z_1|$ ERK has the largest stability regions, followed by EPBM, and lastly EAB. These results suggest that high-order EPBMs will exhibit poor stability properties for purely advective or dispersive equations. However, we can remedy the situation by considering the composite method (\ref{eq:leg_method_composite}). In Figure \ref{fig:stability_iteration} we plot the stability regions for the composite method with $q=4,6,8$ and $\kappa=0,1,2$ where $\kappa$ is the number of times the iterator method is applied. The computational cost of the composite method scales by a factor of $(\kappa+1)$ compared to the propagator, therefore it is critical that the stability regions increases by a similar factor or the composite method will not be as efficient as the propagator alone. Fortunately, by applying even a single iterator the linear stability properties increase drastically, suggesting that composite schemes are suitable for efficiently solving advective or dispersive equations. 

Finally, since the parameter choices and plot axes in Figure \ref{fig:stability_large} are identical to those presented in \cite{buvoli2019esdc} we can also compare the stability regions of EPBMs to those of high-order exponential SDC methods constructed with Chebyshev nodes. Overall ESDC methods have significantly improved stability regions compared to all methods shown in this paper. Moreover, unlike EPBM or EAB, the stability regions of high-order ESDC methods grow with order. 

\section{Numerical experiments and discussion}
\label{sec:numerical_experiments}

We investigate the efficiency of the Legendre EPBMs from Table \ref{tab:leg_methods} by conducting six numerical experiments consisting of four partitioned systems and two unpartitioned systems. The systems all arise from solving partial differential equations using the method of lines. In each case, the reference solutions were computed numerically by averaging the outputs of multiple methods that were run using a small timestep; the averaging prevents any one method from appearing to converge past discretization error.

We divide our numerical experiments and the subsequent discussion into a partitioned and unpartioned section. In both cases we compare Legendre EPBMs against exponential Adams Bashforth (EAB) \cite{beylkin1998ELP, tokman2006efficient}, exponential spectral deferred correction (ESDC) \cite{buvoli2019esdc}, and the fourth-order exponential Runge-Kutta (ERK) methods from \cite{cox2002ETDRK4} and \cite{michels2017stiffly}. Each of these methods pertains to a different method family: EAB methods are linear multistep methods, ESDC are one-step methods, and EPBM are block methods. We selected EAB and ESDC schemes because both families can be constructed at arbitrarily high orders-of-accuracy for both partitioned and unpartitioned problems. This allows a fair comparison across a range of orders. As a reference, we also include fourth-order exponential Runge-Kutta schemes since they have proven very efficient compared to a host of other methods \cite{michels2017stiffly, rainwater2016new, montanelli2016solving}.

We provide the code that was used to run the experiments in \cite{BuvoliEPBMZenodo}. All the timing results presented in this paper were produced on a 12 core, 2.2 Ghz Intel Xeon E5-2650 v4 processor with hyper-threading enabled. Finally, since the equations and spatial discretizations are identical to those presented in  \cite{buvoli2019esdc}, we only list the critical information for each problem below.

	\subsection{Partitioned numerical experiments} 
	
	We consider the following equations that are each equipped with periodic boundary conditions:
		
	\vspace{0.5em}
	\begin{enumerate}[leftmargin=*]
		\item The {\bf Kuramoto-Sivashinsky} (KS) equation from \cite{kuramoto1976persistent, KassamTrefethen05ETDRK4},
\begin{align}
& \frac{\partial u}{\partial t} = -\frac{\partial^2 u}{\partial x^2} -\frac{\partial^4 u}{\partial x^4} - \frac{1}{2} \frac{\partial }{\partial x} \left(u^2\right) \label{eq:kuramoto}, \\[0.5em]
& u(x,t=0) = \cos\left( \tfrac{x}{16} \right) \left( 1 + \sin\left( \tfrac{x}{16} \right) \right), \hspace{1em} x\in[0, 64\pi], \nonumber
\end{align}
integrated to time $t=60$ using a 1024 point Fourier spectral discretization in $x$.

\vspace{0.5em}
\item The {\bf Nikolaevskiy} equation from \cite{simbawa2010nikolaevskiy, grooms2011IMEXETDCOMP},
\begin{align}
& \frac{\partial u}{\partial t} = \alpha \frac{\partial^3u}{\partial x^3} + \beta \frac{\partial^5 u}{\partial x^5} -\frac{\partial^2}{\partial x^2} \left( r - \left( 1 + \frac{\partial^2}{\partial x^2}\right)^2 \right) u -\frac{1}{2}\frac{\partial}{\partial x}\left( u^2 \right) \label{eq:nikolaevskiy}, \\[0.5em]
& u(x,t=0) = \sin(x) + \epsilon \sin \left(\tfrac{x}{25}\right), \hspace{1em} x\in[-75\pi, 75\pi], \nonumber
\end{align}
where $r=1/4$, $\alpha = 2.1$, $\beta = 0.77$, and $\epsilon = 1/10$. The equation is integrated to time $t=50$ using a 4096 point Fourier spectral discretization.

\vspace{0.5em}
\item The {\bf Korteweg-de Vries} (KDV) equation from \cite{zabusky1965interaction}
\begin{align}
& \frac{\partial u}{\partial t} = -\left[ \delta \frac{\partial^3 u}{\partial x^3} + \frac{1}{2}\frac{\partial}{\partial x}(u^2) \right] \\
& u(x,t=0) = \cos(\pi x), \hspace{1em} x \in [0,2], \nonumber
\end{align}
where $\delta = 0.022$. This equation is integrated to time $t=3.6/\pi$ using a 512 point Fourier spectral discretization. 

\vspace{0.5em}
\item The {\bf barotropic quasigeostrophic} equation on a $\beta$-plane with linear Ekman drag and hyperviscous diffusion of momentum from \cite{grooms2011IMEXETDCOMP}:
\begin{align}
\partial_t \nabla^2 \psi &= -\left[\beta \partial_x \psi  + \epsilon \nabla^2 \psi + \nu \nabla^{10} \psi + \mathbf{u} \cdot \nabla (\nabla^2 \psi) \right] \label{eq:quasigeostrophic} \\
\psi(x,y,t=0) &= \frac{1}{8} \exp\left(-8\left(2y^2+x^2/2 - \pi/4 \right)^2 \right), \quad (x,y) \in[-\pi,\pi], \nonumber
\end{align}
 where $\beta = 10$, $\epsilon = 1/100$, $\nu = 10^{-14}$. The equation is integrated to time $t=5$ using a $512\times 512$ point Fourier discretization.
 
\end{enumerate} 
 
Each of these four equations are solved in Fourier space. Since linear derivative operators diagonalize in Fourier space, each initial value problem has the form $\mathbf{y}' = \mathbf{L}\mathbf{y} + N(t,\mathbf{y})$ where $\mathbf{L}$ is a diagonal operator that includes all the discretized linear derivative terms. For example for the KDV equation $\mathbf{L} = \text{diag}(i\delta \mathbf{k}^3)$ where $\mathbf{k}$ is a vector of Fourier wave numbers. Computing matrix-vector products with $\varphi$-functions of $\mathbf{L}$ now amounts to multiplication with a diagonal matrix, or equivalently element-wise vector multiplication. Moreover, since our experiments only consider constant stepsizes, we can efficiently precompute and store the diagonal entries of $\varphi_j(\mathbf{L})$ as a vector. For small wave numbers with magnitude less than one we compute the diagonal entries of $\varphi_j(\mathbf{L})$ using the contour integral method \cite{KassamTrefethen05ETDRK4}, while for larger wave numbers we simply  use the recurrence relation ${\varphi_{j+1}(z) = z^{-1}(\varphi_j(z) - 1/j!), ~ \varphi_0(z) = e^z}$. Lastly, for antialiasing we apply the standard two-thirds rule to zeros out the top third of the spectrum.

One of the key features of Legendre EPBMs is their inherent parallelism that allows each of the method's outputs to be computed simultaneously. In order to investigate the potential advantages of this approach we implemented all the partitioned integrators in Fortran and used OpenMP to parallelize the EPBM timestep computation. The code was therefore run using a single thread for serial methods, and using $q$ threads for a Legendre method with $q$ nodes. 

Lastly, in Appendix \ref{ap:complex_nodes_epbms} we also briefly compare EPBMs with imaginary equispaced nodes to the EPBMs with Legendre nodes on the Kuramoto equation.

\begin{figure}%
\centering
\includegraphics[trim={4cm 8.5cm 4cm 1cm},clip,width=1\linewidth]{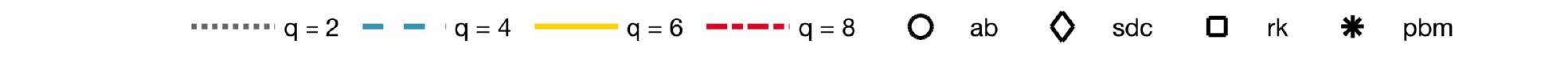}
\includegraphics[trim={1.1cm 0.5cm 4.25cm 9.4cm},clip,width=1\linewidth]{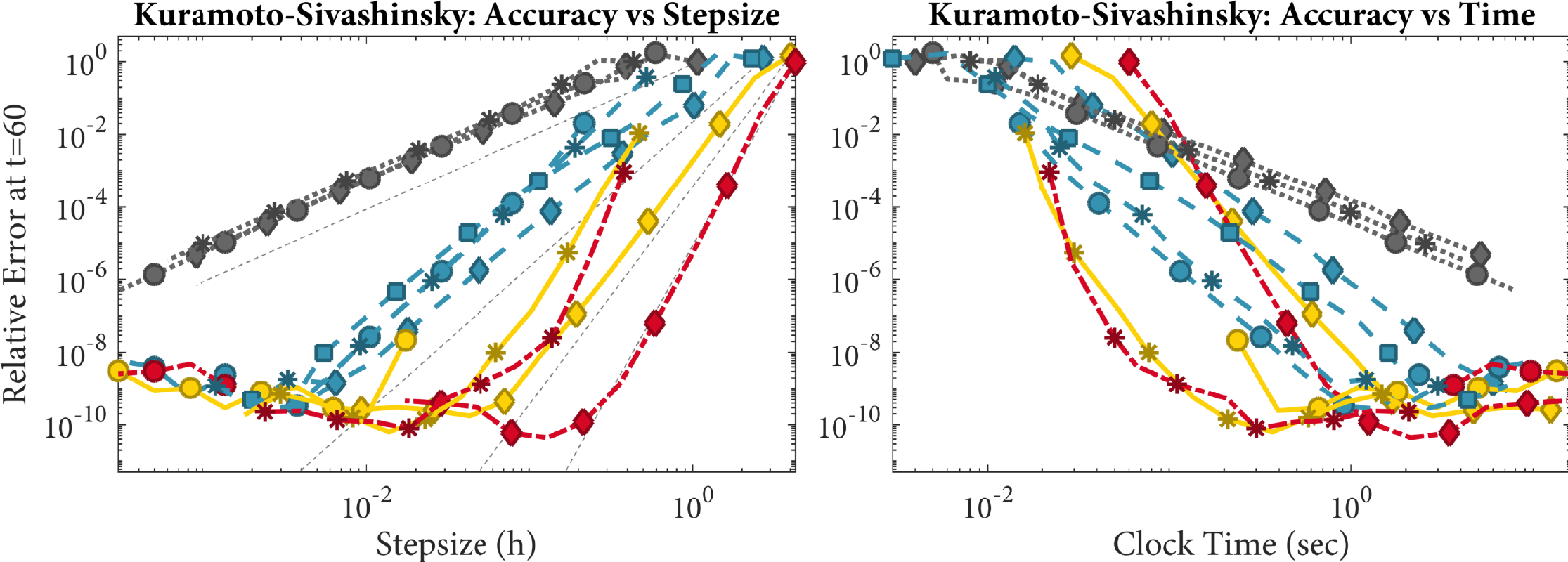}
\includegraphics[trim={1.1cm 0.5cm 4.25cm 9.4cm},clip,width=1\linewidth]{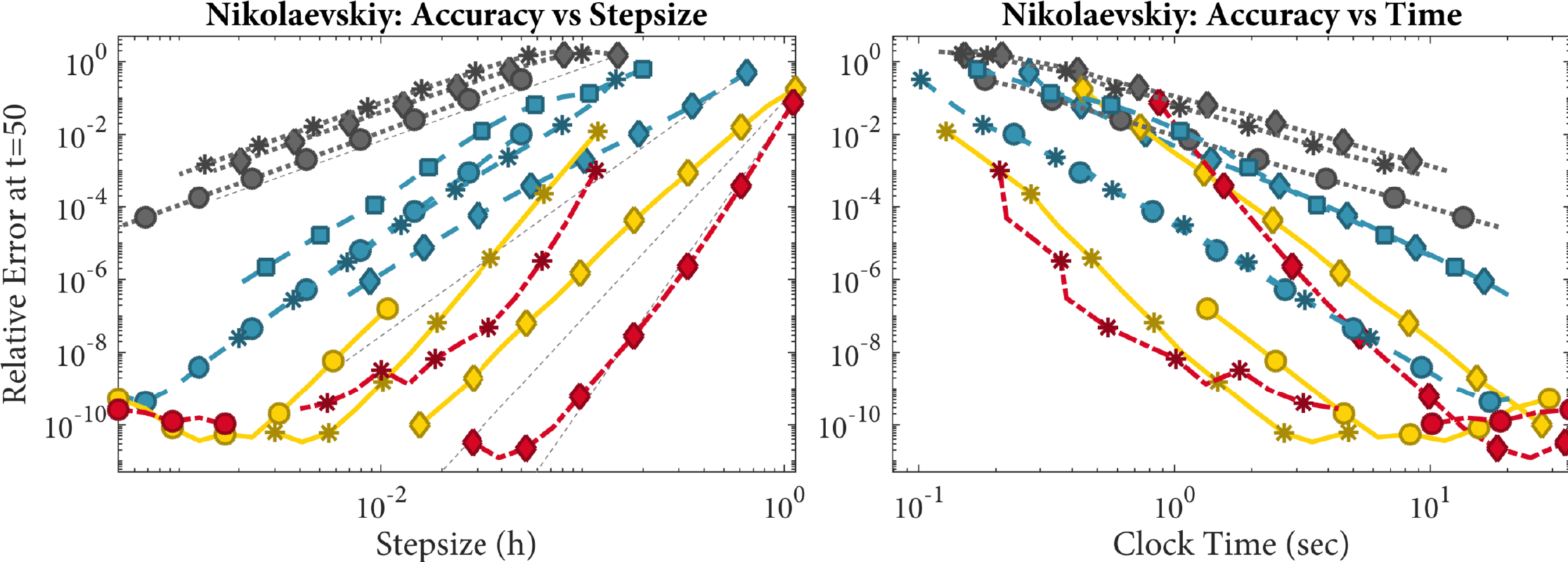}
\includegraphics[trim={1.1cm 0.5cm 4.25cm 9.4cm},clip,width=1\linewidth]{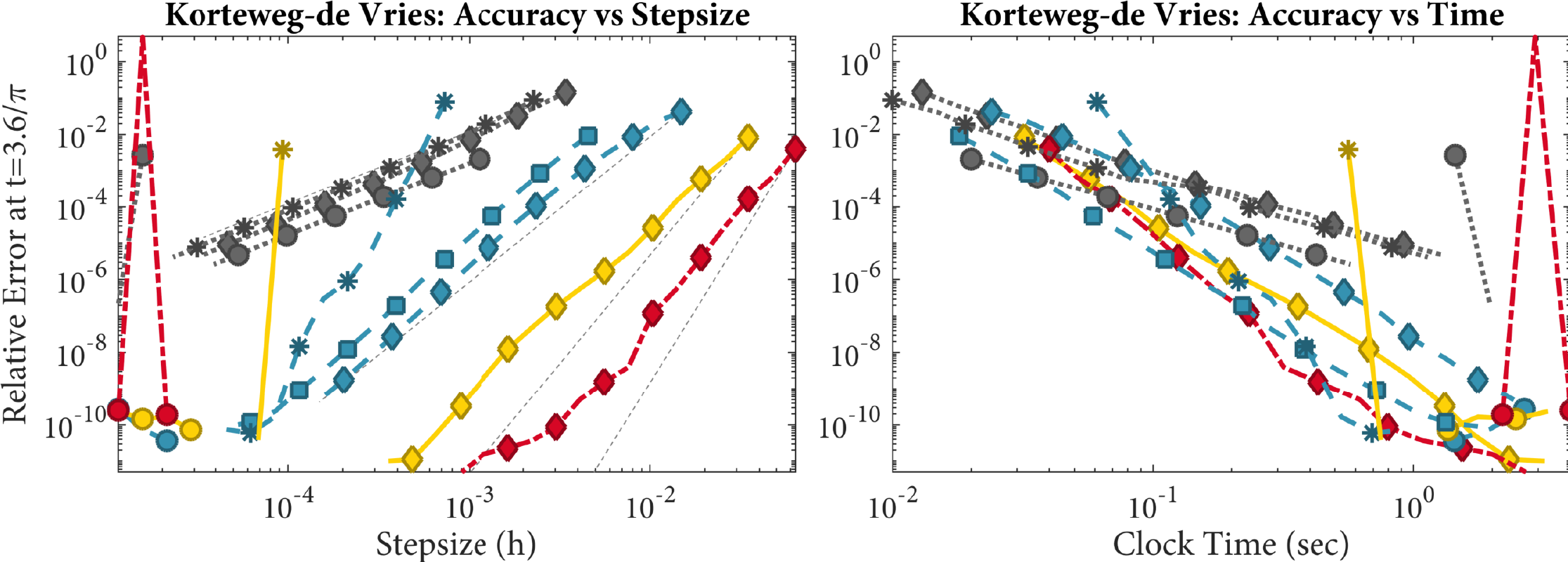}
\includegraphics[trim={1.1cm 0.5cm 4.25cm 9.4cm},clip,width=1\linewidth]{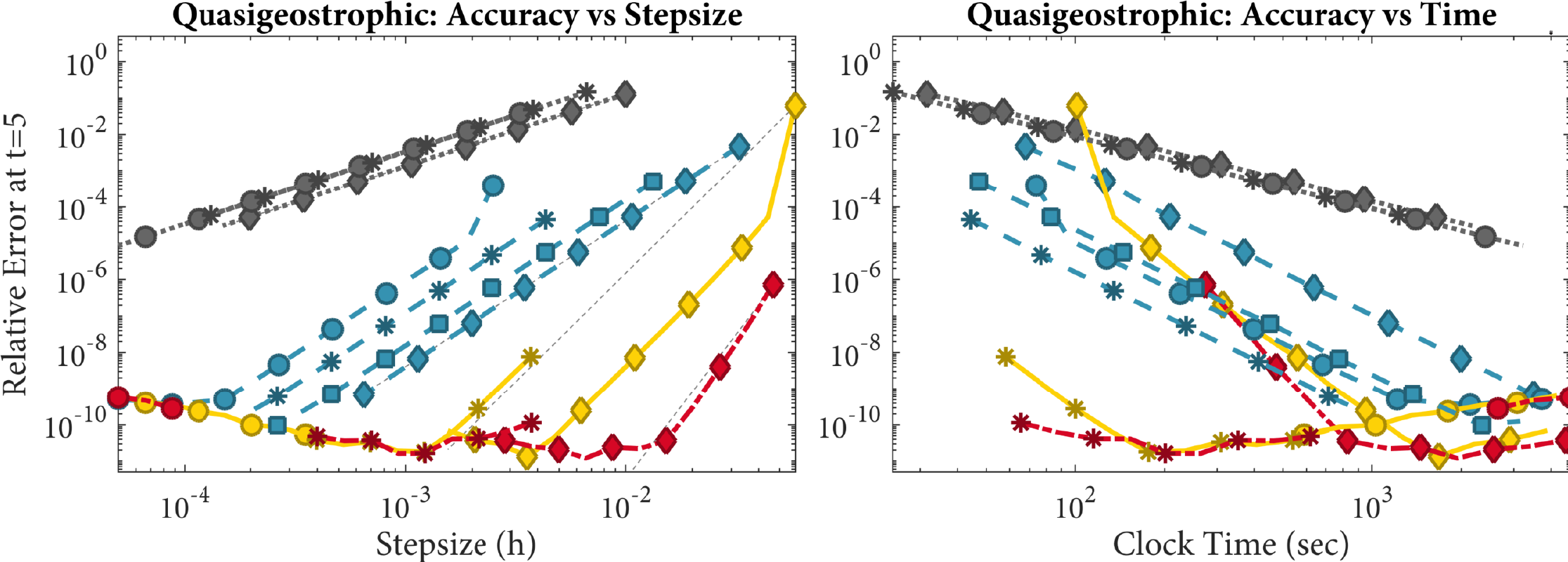}

\vspace{-1em}
\caption{Plots of accuracy versus stepsize and accuracy versus clock-time for the four partitioned equations. All the methods shown in this figure are exponential integrators. Color represents order-of-accuracy and the markers denote the method type. In the legend pbm stands for polynomial block methods, sdc for spectral deferred correction, rk for Runge-Kutta, and ab for Adams Bashforth. The thin, dashed gray lines in the accuracy versus stepsize plots correspond to second, fourth, sixth, and eight order convergence.}
\label{fig:results_page1}
\end{figure}
\begin{figure}%
\centering
\includegraphics[trim={13cm 0.5cm 4.25cm 9.4cm},clip,width=0.49\linewidth]{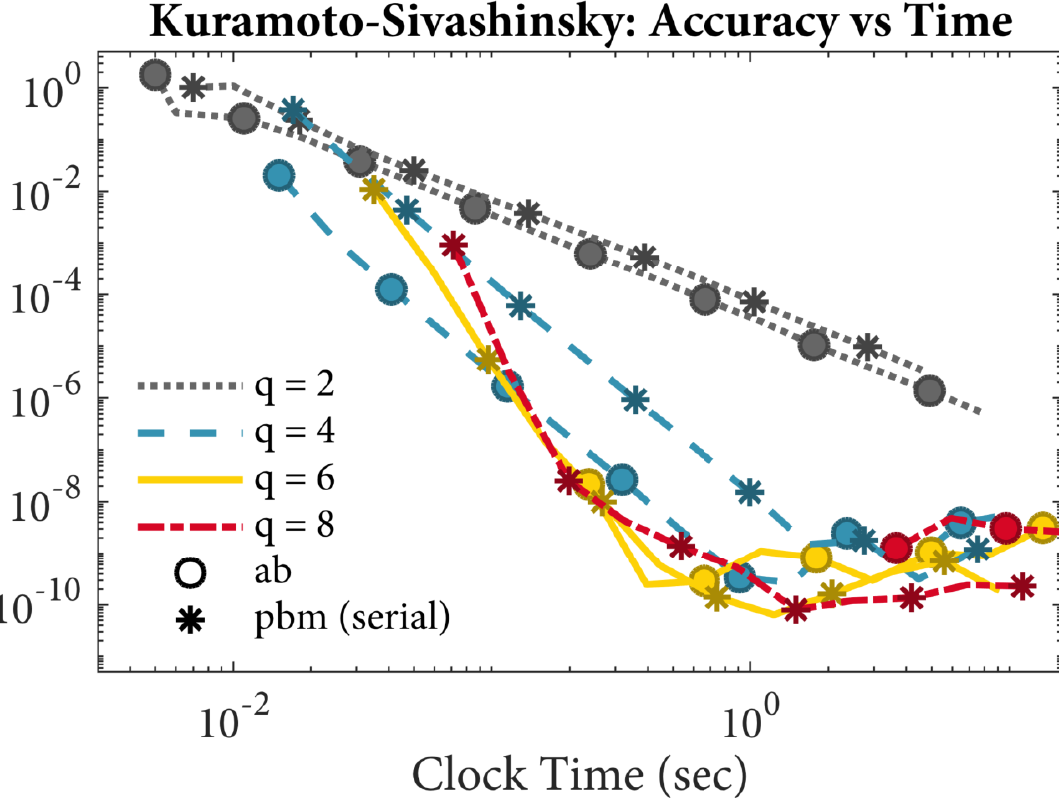}
\includegraphics[trim={13cm 0.5cm 4.25cm 9.4cm},clip,width=0.49\linewidth]{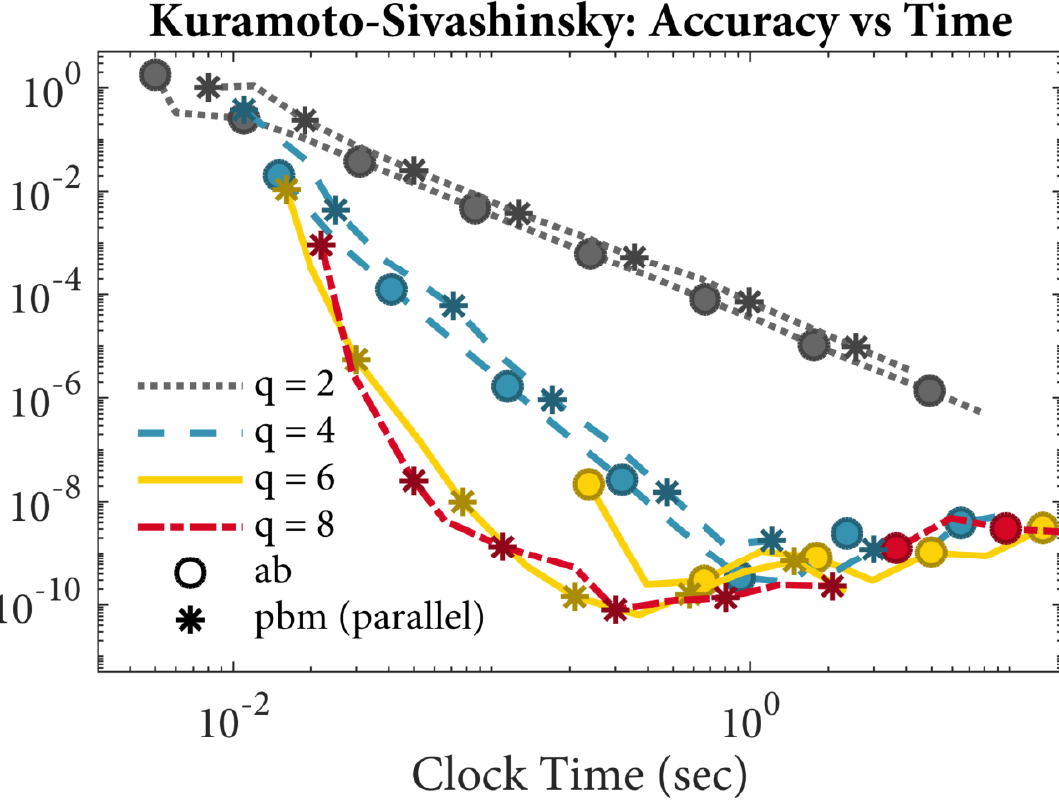}
\caption{Quantifying the effects of parallelism for EPBMs on the KS equation. The two plots show the accuracy versus clock time for EAB and EPBMs that are implemented in serial (left) and in parallel (right). The right plot is identical to the one shown in Figure \ref{fig:results_page1} except that ESDC and ERK methods have been removed.}
\label{fig:kuramoto_serial_vs_parallel}
\end{figure}

\begin{figure}%
\centering
\includegraphics[trim={4cm 8.5cm 4cm 1cm},clip,width=1\linewidth]{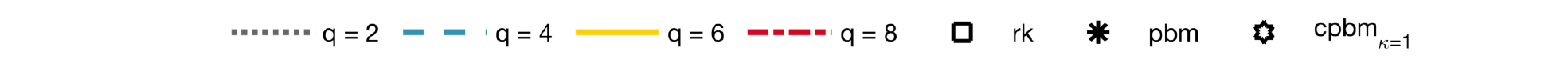}
\includegraphics[trim={1.1cm 0.5cm 4.25cm 9.4cm},clip,width=1\linewidth]{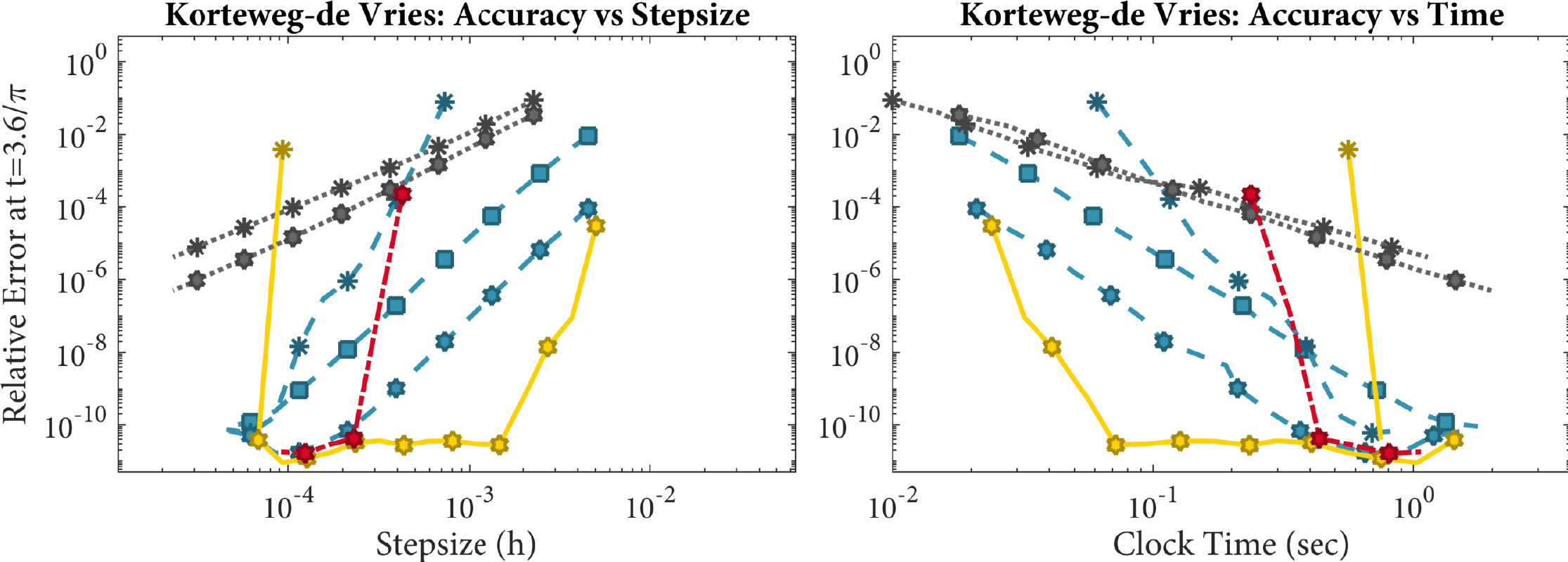}
\caption{Composite methods for solving KDV. We compare Legendre EPBMs to the composite Legendre methods (\ref{eq:leg_method_composite}) with $k=1$. We also include ERK4 which was one of the most efficient methods for solving KDV in Figure \ref{fig:results_page1}. As suggested by the linear stability results, composing a PBM with a single iterator leads to a composite method with drastically improved stability on dispersive equations.}
\label{fig:results_epbm_m}
\end{figure}

\begin{figure}%
	\centering
	\includegraphics[trim={4cm 8.5cm 4cm 1cm},clip,width=1\linewidth]{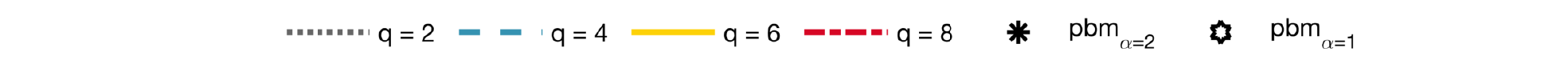}
	\includegraphics[trim={1.1cm 0.5cm 4.25cm 9.4cm},clip,width=1\linewidth]{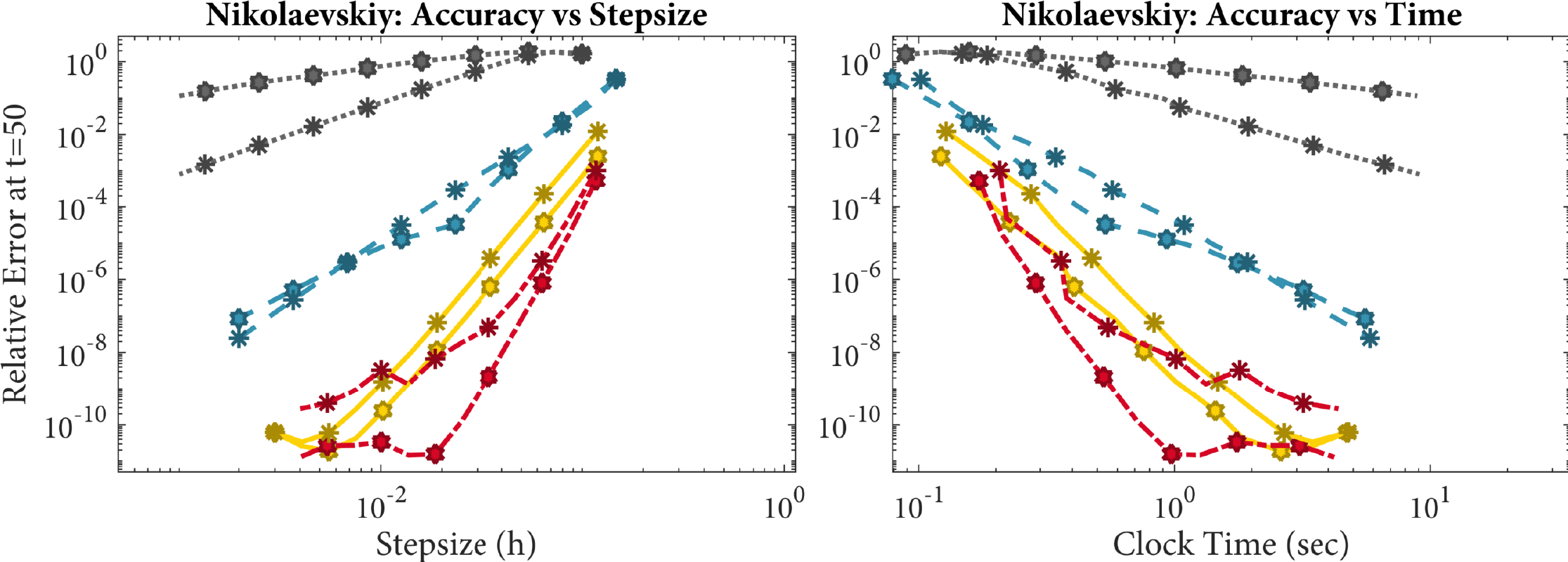}
	\caption{Quantifying the effects of the extrapolation parameter $\alpha$.  Accuracy versus stepsize and accuracy versus clock time for the Nikolaevskiy equation using Legendre EPBMs with $\alpha=2$ and $\alpha=1$. High-order methods benefit from a smaller $\alpha$ to avoid problems with round-off errors, however low order methods have better error constants with larger $\alpha$.}
	\label{fig:results_epbm_alpha}
\end{figure}

\subsection{Discussion of partitioned numerical experiments}		
			
	In Figure \ref{fig:results_page1} we present plots of relative error versus stepsize, and relative error versus computational time for each of the four partitioned problems. In Figure \ref{fig:kuramoto_serial_vs_parallel} we also investigate the efficiency of EPBMs when they are run in serial on a single thread by repeating our Kuramoto experiment. For the dispersive KDV equation we conduct a separate test in Figure \ref{fig:results_epbm_m} using the composite Legendre methods (\ref{eq:leg_method_composite}) in order to validate the results of the linear stability analysis. Finally, in Figure \ref{fig:results_epbm_alpha} we investigate the effects of varying $\alpha$ on the Nikolaevskiy equation where high-order EPBMs showed order reduction.
		
		Stability, accuracy, and computational cost per timestep are the three factors that determine a method's overall efficiency. An ideal method should have a large stability region that allows for coarse stepsizes, a small error constant to outcompete other methods of the same order, and a low computational cost per timestep. Each of the three method families presented in our numerical experiments has different strengths and weaknesses. EAB methods have the lowest cost per timestep, but they suffer from limited stability and poor error constants at high orders. EPBMs have improved error constants and stability compared to EAB and nearly the same timestep cost when implemented in parallel. ESDC methods have excellent stability and accuracy compared to other methods, however their computational cost is significantly larger. 
		
	The stability region and error-constant are fixed properties of a method, however the cost per timestep is highly problem-dependent. For the partitioned experiments we precompute the matrix exponential functions, and the cost of computing a matrix-vector product $\varphi_j(\mathbf{L})\mathbf{x}$ is equivalent to that of an element-wise vector multiplication between $\mathbf{x}$ and a precomputed vector. Therefore, the majority of the computational work for the partitioned experiments is due to nonlinear function evaluations, and the cost of a method is determined by the number of function evaluations required per timestep.

	EAB methods only require a single nonlinear function evaluation per timestep and serve as a benchmark for optimal compute cost. In contrast, Legendre EPBMs require one function evaluation per node. However, since the outputs are independent, each function evaluation can be computed in parallel. If we neglect parallel overhead, then EPBMs have the same cost per timestep as an EAB method. On the other hand an ESDC method with $q$ nodes requires $q(q-1)$ serial function evaluations per timestep, making the cost significantly higher than that of the other methods. Finally the fourth-order Runge-Kutta method requires four serial function evaluations per timestep.
	
	Amongst the second-order methods, EAB is consistently most efficient due to it's superior accuracy compared to the other integrators. At fourth order, the results were more variable. On the Kuramoto and Nikolaevskiy equation, EAB was the most efficient at fine error tolerances, while at coarse tolerances EPBM was more efficient since EAB was no longer stable. For KDV the fourth order ETDRK4 was the most efficient except near discretization error where EPBM is stable. At orders six and eight, the Legendre EPBMs were able to consistently achieve the best accuracy in the smallest amount of computational time. This was due to their improved accuracy and use of parallelism. For one-dimensional problems with dissipation, EPBMs were significantly more efficient than existing methods. However, when solving Kuramoto and Nikolaevskiy, EPBM8 showed order reduction at small timesteps. This phenomenon is due to rounding error and is present in other high-order polynomial integrators \cite{buvoli2019constructing, buvoli2018polynomial}. As is the case for classical integrators, the issue can be resolved by simply selecting a smaller extrapolation factor. As shown in Figure \ref{fig:results_epbm_alpha}, reducing $\alpha$ to one, makes EPBM8 the most efficient method for obtaining highly accurate solutions.
	
	When EPBMs were run in serial on the Kuramoto equation (Figure \ref{fig:kuramoto_serial_vs_parallel}), their computational advantages decreased significantly. Fourth-order EPBM is now less efficient than the fourth-order EAB, while the sixth and eight order EPBMs only marginally outperformed EAB schemes at the lowest error tolerances. This behavior is expected because a serial implementation of a Legendre method with $q$ nodes requires $q-1$ function evaluations per timestep, compared to the single function evaluation for an EAB scheme. In short, Figure \ref{fig:kuramoto_serial_vs_parallel} shows that the large computational gains for EPBM schemes are only possible if the methods can be parallelized efficiently.
	
	Finally, as predicted by linear stability analysis, EAB methods and Legendre EPBMs demonstrated poor stability on the dissipative KDV equation and only converged at smaller time-steps. However, the composite Legendre method with $k=1$ shown in Figure \ref{fig:results_epbm_m} was able to retain stability across a much wider range of timesteps. Moreover, though its cost per timestep is twice that of a Legendre EPBM, the improved stability and accuracy made both the fourth-order and sixth order EBPM significantly more efficient than ETDRK4. 
	
	When comparing methods across different orders-of-accuracy we see that high-order EPBMs are frequently able to obtain extremely accurate solutions in the same amount of time that it takes a second-order integrator to achieve only a few digits of accuracy. Moreover, if we fix the timestep, we see that all stable EPBMs run in nearly the same computational time. The equations presented here are simple and the domains are conveniently periodic so that we may diagonalize the linear operators. Nevertheless our result suggest that for similar equations, or more broadly for problems where efficient parallelism is possible, low-order EPBMs should only be considered if one wants the solution in the fastest possible time and requires very little accuracy; if more accuracy is required then, barring other computational restrictions such as limited memory, a high-order EPBM can obtain a highly accurate solution at little to no extra cost.
	
\subsection {Unpartitioned initial value problems} 
	
We consider the two dimensional ADR equation with homogeneous Neumann boundary conditions from \cite{rainwater2016new},
\begin{align}
&u_t = \epsilon\left(u_{xx} + u_{yy}\right) + \delta \left( u_{x} + u_{y}\right) + \gamma u(u-1/2)(1-u), \\
&u(x,t=0) = 256(xy(1-x)(1-y))^2 + 0.3,	\nonumber \\ \nonumber
&x,y \in [0, 1].
\end{align}
The equation is integrated to time $t=1/100$ and the spatial discretization employs standard second-order accurate finite differences on a 200$\times$200 point grid. We solve this equation using the two parameter sets described in Table \ref{tab:unpartitioned_spectral_radius} that lead to either a stiff linearity or a stiff nonlinearity.

\begin{table}[h]
	
	\begin{center}
	\renewcommand{\arraystretch}{1.25}
	\begin{tabular}{r|lll|ll}
		& $\epsilon$ & $\delta$ & $\gamma$ & $\rho(\mathbf{L})$ & $\rho\left( \frac{\partial N}{\partial y}\right)$ \\ \hline
		stiff linearity & $1/100$ & $-10$ & $100$ & $4288$ & $167$\\
		stiff nonlinearity & $1/10000$ & $-1 / 10$ & $1000$ & $43$ & $1670$ 
	\end{tabular}
	\end{center}
	\vspace{0.5em}
	
	\caption{Parameters for the 2D ADR equation. We also show the spectral radius $\rho$ for the linear and nonlinear terms.}
	\label{tab:unpartitioned_spectral_radius}
	
\end{table}

\begin{figure}[!htb]
\centering

	\begin{footnotesize}
		{\bf  (a)}  ADR - stiff linearity
	\end{footnotesize}
	\vspace{1em}
	
	\includegraphics[trim={1cm 10.75cm 2.75cm 0.5cm},clip,width=1\linewidth]{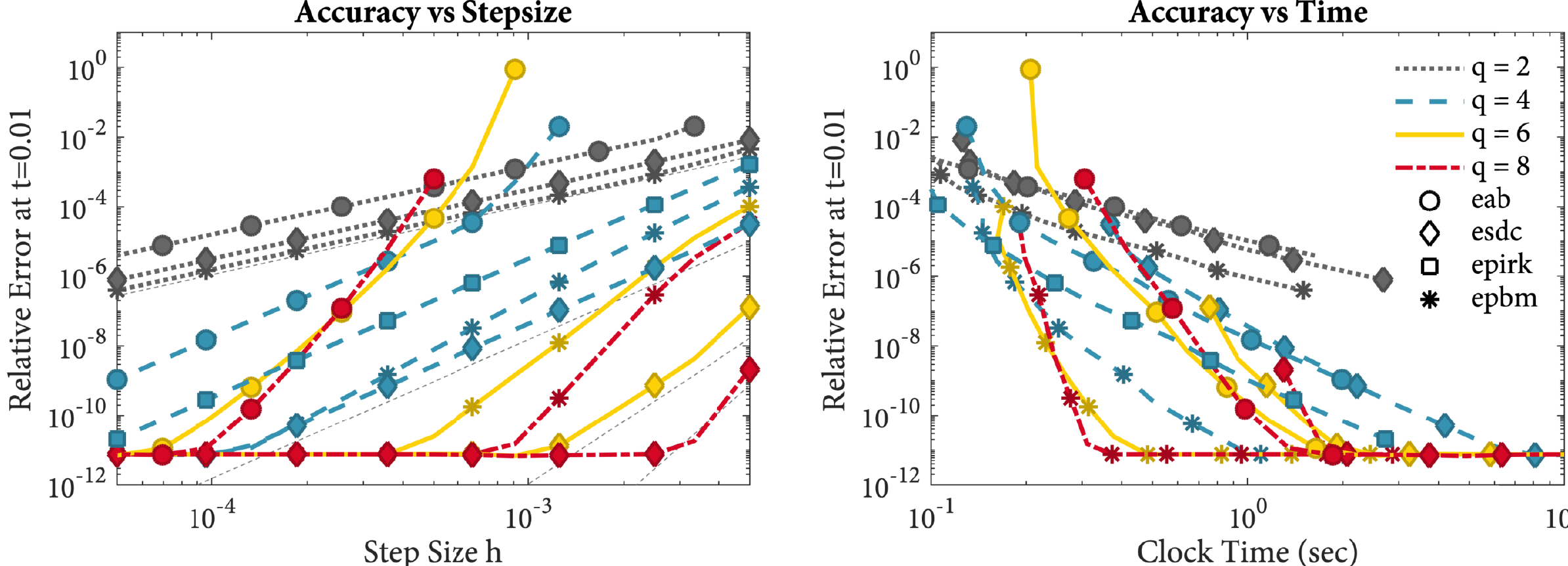}

	\begin{footnotesize}
	{\bf (b) }  ADR - stiff nonlinearity	
	\end{footnotesize}
	\vspace{1em}

	\includegraphics[trim={1cm 10.75cm 2.75cm 0.5cm},clip,width=1\linewidth]{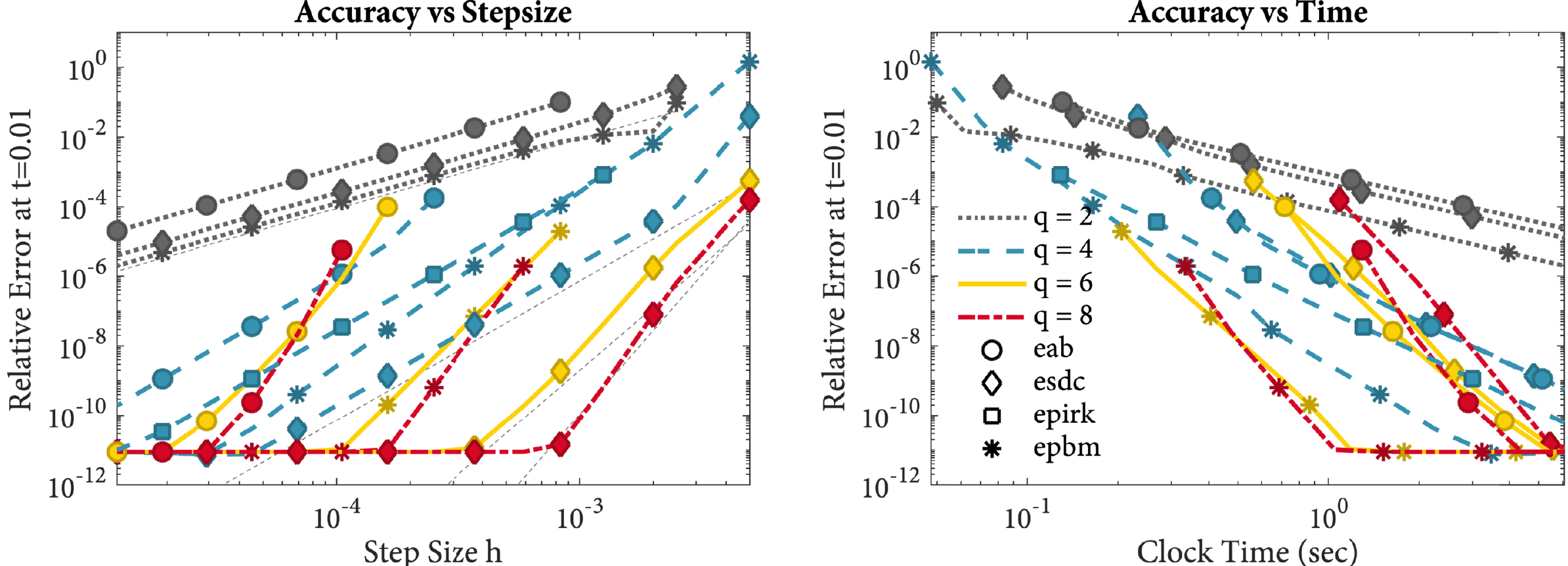}

\caption{(a) Diagrams for the ADR equation with a stiff linearity. (b) Diagrams for the ADR equation with a stiff nonlinearity. The dashed lines in the accuracy diagrams correspond to second, fourth, sixth and eighth order convergence.
}
\label{fig:results_page3}
\end{figure}

To compute matrix-vector products with $\varphi$-functions we apply the Krylov-based KIOPs algorithm \cite{GAUDREAULT2018236} using the exact Jacobian. In Figure \ref{fig:results_page3} we show results for Legendre EPBMs of orders two, four, six and eight compared against unpartitioned Adams-Bashforth methods \cite{tokman2006efficient}, unpartitioned exponential spectral deferred correction methods \cite{buvoli2019esdc}, and the fourth-order unpartitioned exponential Runge-Kutta method EPIRK43s \cite{michels2017stiffly}.

\subsubsection{Discussion of unpartitioned numerical experiments}

	All unpartitioned exponential integrators require matrix-vector products with $\varphi$-functions of the Jacobian. Since the Jacobian varies in time, we are no longer able to store the exponential matrix functions at the first time step. Moreover, the Jacobian matrices are non-diagonal, therefore the majority  of the computational effort at each time step is now due to the vector products with $\varphi$-functions.
				
	 Krylov methods like KIOPS \cite{GAUDREAULT2018236} and PHIPM \cite{NiesenWright2012Krylov, NiesenWright2011Krylov} approximate matrix-vector products of $\varphi$-functions within a Krylov subspace. %
	 The operation to build a Krylov space and approximate a $\varphi$-function is called a {\em projection} and its cost will depend on the spectrum of $\mathbf{A}$. The cost of a time-integrator is therefore dependent on the total number of projections per time step. In general, every $\varphi$-function in an arbitrary linear combination must be computed using a separate projection. However, expressions of the form 
		\begin{align}
			\varphi_0(\tau \mathbf{A})\mathbf{x}_0 + \tau \varphi_1(\tau \mathbf{A})\mathbf{x}_1 + \tau^2 \varphi_2(\tau \mathbf{A})\mathbf{x}_2 + \ldots + \tau^p \varphi_p(\tau \mathbf{A})\mathbf{x}_p
			\label{eq:phi_projection}
		\end{align}
 		can be evaluated at multiple $\tau$ values in a single projection. Furthermore, this computation is made cheaper if $\tau$ is small. This fact can be used to construct highly efficient integrators such as EPIRK methods \cite{tokman2006efficient, tokman2011new, rainwater2016new, michels2017stiffly} whose coefficients are specifically chosen to minimize the total number of projections required at each timestep. 
 		
 		It is important to state that Krylov methods are not the only way to compute the expression (\ref{eq:phi_projection}). Leja interpolation \cite{caliari2016leja} and scaling-and-squaring based algorithms \cite{higham2020catalogue} are two such alternatives; a detailed comparison of these methods is presented in \cite{caliari2014comparison}, and \cite{higham2020catalogue} contains a survey of existing codes.
 		
 		In our experiment we used KIOPs to compute the exponential terms.  Exponential Adams-Bashforth methods once again serve as a benchmark for optimal cost per timestep as they only require one projection where $\tau = h$. In contrast, EPIRK43s requires two serial projections: the first has $\tau = h$ and the second has $\tau = h/9$. For polynomial integrators, all \ppes{} can be written in the form (\ref{eq:phi_projection}) and therefore each output of any EPBM requires one projection to compute. However, since Legendre EPBMs are constructed using a single Adams \ppe{}, all the outputs can be computed in one projection where $\tau = 2h$.  Finally, ESDC methods are the most expensive requiring one projection per correction iteration, meaning that a method with $q$ nodes requires $q$ projections to $\tau=h$. Under these conditions, the cost per timestep of an unpartitioned Legendre EPBM can be up to roughly twice as expensive as EPIRK43s and EAB methods, and a factor of $q/2$ times less expensive than an ESDC method with $q$ nodes. Finally, we also note that for unpartitioned Legendre EPBMs it is possible to parallelize the right-hand-side evaluations. However, for our test problems the cost of right-hand-side evaluations of $R(\mathbf{y})$ is negligible compared to that of the Krylov projections, hence we always evaluate the right-hand-side in serial. 	
 	
 		Similar to our partitioned experiments we see that high-order ESDC methods have the best stability and error constants, however their high computational cost ultimately renders them less efficient. On the other hand, EPBMs possess improved stability and accuracy compared to EAB schemes and were overall the most efficient methods at each order-of-accuracy. Amongst the non-polynomial integrators, the fourth-order EPIRK method outperformed all other EAB and ESDC integrators except at very high precision where 8th-order EAB becomes more efficient.
 		
 		Unlike for the partitioned experiments, we see that higher-order unpartitioned EPBMs are only more efficient than low-order EPBMs if high-accuracy is required. In particular it would be wasteful to always use a high-order integrator since it will require a substantially finer stepsize (and therefore more computational work) than a low-order method to remain stable. 
 		 		 			
 		Finally we note that it is trivial to construct methods that only require a projection with $\tau=h$ per timestep. For example, we obtain this property by choosing Chebyshev nodes and PMFC$_1$ with endpoints $b_j = 1$. It is also simple to apply parallelism to speed up the computation of methods that use different endpoints $b_j$ or polynomials $L^{[j]}_R(\tau)$ for computing each output. Though such methods will require $q$ projections per timestep, each projection can now be computed in parallel to nullify the additional cost. %

\section{Summary and Conclusions}

We presented an extension of the polynomial framework that incorporates exponential integration. To achieve this we combined the classical \ODEp{} with the  integrating factor method to introduce the \ppes{} that form the basis of the new methods. By utilizing \ppes{} it is possible to construct many families of new exponential integrators including those with desirable properties such as parallelism and high-order. 

After introducing the framework we demonstrated its potential by presenting several general construction strategies for EPBMs and by deriving a new class of parallel EPBMs that use Legendre points. Our numerical experiments demonstrate the potential of these new parallel methods for both partitioned and unpartitioned problems. Based on our results it appears that Legendre EPBMs  can provide significant computational savings compared to current state-of-the-art methods if they can be efficiently parallelized.

The generality of the exponential polynomial framework creates many opportunities for additional developments that we plan to address in future work. In particular we will investigate the construction of exponential polynomial general linear methods that can exploit parallelism on a larger scale. %

\section{Acknowledgments}

I would like to thank Randall LeVeque for encouraging me to develop polynomial integrators during my graduate studies and for his comments and guidance through the course of this project. I would also like to thank Mayya Tokman for her helpful comments on drafts of this work.

\appendix

\section{Coefficients for Adams \ppe}
\label{ap:adams_ppe_coefficients}

Here we describe how to rewrite the \Ape{} from Definition \ref{def:adams_phi_expansion} in terms of the values in an exponential \ODEd{}. In short, one must rewrite $L_y(b)$ and $L_N^{(k)}(b)$ in terms of the values and derivatives that were used to construct the Lagrange interpolating polynomials $L_y(\tau)$ and $L_N(\tau)$. In other words, we must find the coefficients $a_j$ and $c_{\nu,j}$ so that
		\begin{align}
			L_y(b) = \sum_{j=1}^w a_j \mathbf{y}_j \quad \text{and} \quad L^{(\nu)}_N(b) = \sum_{j=1}^w c_{\nu,j} r\mathbf{N}_j.
		\end{align}
		
		If the solution value $\mathbf{y}_k$ is not used to form $L_y(\tau)$, then $a_k$ is zero. All the non-zero $a_j$ are finite difference weights for computing the zeroth derivative at $\tau=b$, using data at the nodes of the Lagrange polynomial  $L_y(\tau)$.
		
		Similarly, if the derivative component value $r\mathbf{N}_k$ is not used to form $L_N(\tau)$, then $c_{\nu, k}=0$, and all the non-zero $c_{\nu,j}$ are finite difference weights for computing the $\nu$th derivative at $\tau=b$, using data at the nodes of the Lagrange polynomials $L_N(\tau)$. 	
		
		The direct procedure for computing the finite difference weights is as follows. Suppose that $L(x)$ is a Lagrange polynomial with nodes $x_j$ and data $\mathbf{g}_j$, for ${j=1, \ldots, l}$. The finite difference weights $w_{\nu,j}$
			for computing the $\nu$th derivative of $L(x)$ at $x=x_0$ (i.e. $L^{\nu}(x_0) = \sum_{j=1}^l w_{\nu,j} \mathbf{g}_j$) are $d_{\nu, j} = \nu! \mathbf{V}^{-1}_{\nu + 1, j}$ where $\mathbf{V}$ is the $l\times l$ Vandermonde matrix with entries $\mathbf{V}_{i,j} = (x_i - x_0)^{j-1}$. This direct procedure can be used to obtain the weights, however, since Vandermonde matrices can be ill-conditioned it is advantageous to use the fast, stable algorithm developed by Fornberg for computing finite difference weights \cite{fornberg1988, fornberg1998classroom}.

\section{Method Coefficients}
\label{ap:leg_method_coeff}

The general form for a Legendre method from Table \ref{tab:leg_methods} with $q$ nodes is
\begin{align}
	y^{[n+1]}_j = \varphi_0(r \eta_j \mathbf{L}) y^{[n]}_1 + r \sum_{k=1}^{q-1} \eta_j^{k} \varphi_{k}(r \eta_j \mathbf{L}) \mathbf{v}_j, && j = 1, \ldots q,
\end{align}
where $\eta_j = z_j + \alpha + 1$ and $\mathbf{N}_j = N(z_j, \mathbf{y}(z_j))$. The vectors $\mathbf{v}_j$ for can be derived by computing the derivatives of the Lagrange polynomial 
\begin{align}
	L(\tau) = \sum^{q-1}_{j=1} \ell_j(\tau) \mathbf{N}_{j}, && \ell_j(\tau) &= \prod^{q-1}_{\substack{l=1\\l \neq j}} \frac{\tau-x_l}{x_j - x_l},
\end{align}
at the point $\tau = -1$, where $x_j$ is the $j$th zero of the $q$th Legendre polynomial $P_{q}$(x). Below we provide several formulas for the $\mathbf{v}_j$ if $q=2,3,4.$

\vspace{1em}
\begin{center}
	\begin{minipage}[t]{0.75\textwidth}
		{\bf Legendre Method:} $q=2$, $z_j = \left\{ -1,  0 \right\}$ 
		\begin{align*}
		\mathbf{v}_1 = \mathbf{N}_2
		\end{align*}
		
		{\bf Legendre Method:} $q=3$, $z_j = \left\{ -1,  -\sqrt{\frac{1}{3}}, \sqrt{\frac{1}{3}} \right\}$
		\begin{align*}
		\mathbf{v}_1 &= \frac{(1+\sqrt{3})\mathbf{N}_2}{2}  + \frac{(1-\sqrt{3}) \mathbf{N}_3}{2}  \\
		\mathbf{v}_2 &= \frac{\sqrt{3} \mathbf{N}_3}{2}-\frac{\sqrt{3} \mathbf{N}_2}{2}
		\end{align*}
		
		{\bf Legendre Method:} $q=4$, $\{z_j\} = \left\{ -1,  -\sqrt{\frac{3}{5}}, 0, \sqrt{\frac{3}{5}} \right\}$ 
		\begin{align*}
		\mathbf{v}_1 &= \frac{(5+\sqrt{15})\mathbf{N}_2}{6} -\frac{2 \mathbf{N}_3}{3}-\frac{(\sqrt{15}-5)\mathbf{N}_4}{6} \\
		\mathbf{v}_2 &= \frac{(-10-\sqrt{15})\mathbf{N}_2}{6} +\frac{10 \mathbf{N}_2}{3}+\frac{(\sqrt{15}-10)\mathbf{N}_4}{6} \\
		\mathbf{v}_2 &= \frac{5 \mathbf{N}_2}{3}-\frac{10 \mathbf{N}_3}{3}+\frac{5 \mathbf{N}_4}{3}	
		\end{align*}
	\end{minipage}
\end{center}
\vspace{1em}

\section{Numerical experiment for EPBMs with complex nodes}
\label{ap:complex_nodes_epbms}

\begin{figure}[h]
	\includegraphics[trim={4cm 8.5cm 4cm 1cm},clip,width=1\linewidth]{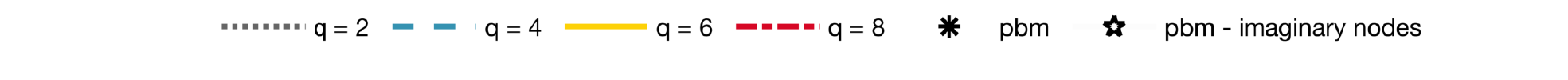}
	\includegraphics[trim={1.1cm 0.5cm 4.25cm 9.4cm},clip,width=1\linewidth]{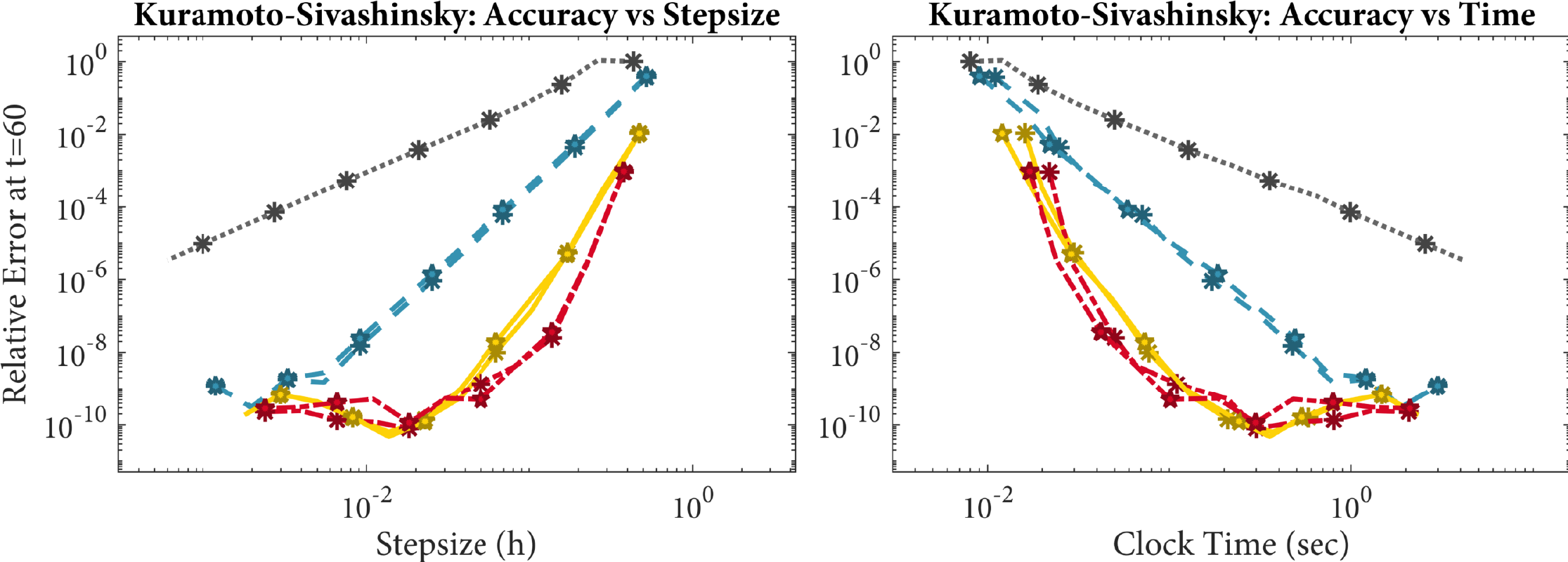}
	\caption{Comparison between methods with real-valued Legendre nodes and imaginary equispaced nodes on the KS equation.}
	\label{fig:kuramoto_iequi}
\end{figure}

In Figure \ref{fig:kuramoto_iequi}, we briefly demonstrate the ability to use EPBMs with complex nodes by solving the Kuramoto equation (\ref{eq:kuramoto}) using the EPBMs with Legendre nodes from Table \ref{tab:leg_methods} and the EPBMs with imaginary nodes labeled as ``method A'' from Table \ref{tab:iequi_methods}. Both methods are run with $\alpha=2$. Though they offer no computational advantage, the methods with imaginary nodes retain stability and have similar error properties to the real-valued methods on this equation. 

\vspace{0.5em}

\bibliographystyle{siamplain}
\bibliography{references_nourl}

\end{document}